\documentclass[11pt]{tran-l}

\usepackage{fullpage}
\usepackage[centertags]{amsmath}
\usepackage{amsfonts}
\usepackage{amssymb}
\usepackage{amsthm}
\usepackage[all]{xy}

 \newtheorem{thm}{Theorem}[section]
 \newtheorem{cor}[thm]{Corollary}
 \newtheorem{lem}[thm]{Lemma}
 \newtheorem{prop}[thm]{Proposition}
 \theoremstyle{definition}
 \newtheorem{defn}[thm]{Definition}
 \theoremstyle{remark}
 \newtheorem{rmk}[thm]{Remark}
 \numberwithin{equation}{section}
 \newcommand{\To}{\longrightarrow}

 \newcommand{\cD}{\mathcal{D}}
 \newcommand{\cR}{\mathcal{R}}
 \newcommand{\cX}{\mathcal{X}}
 \newcommand{\cA}{\mathcal{A}}
 \newcommand{\cH}{\mathcal{H}}
 
 \newcommand{\cL}{\mathcal{L}}
 \newcommand{\cF}{\mathcal{F}}
 \newcommand{\cO}{\mathcal{O}}

 \newcommand{\bR}{\mathbb{R}}
 \newcommand{\bP}{\mathbb{P}}
 \newcommand{\bC}{\mathbb{C}}

 \newcommand{\bZ}{\mathbb{Z}}
 \newcommand{\bQ}{\mathbb{Q}}
 
 \newcommand{\bN}{\mathbb{N}}

 \newcommand{\Hcb}{H^1_c(\Gamma_{0}(Np),\mathcal{D})}
 \newcommand{\Hc}{H^1_c(\Gamma_{0}(Np),\mathcal{D}_N)}
 \newcommand{\HB}{H^1_c(\Gamma_{0}(Np),\mathcal{D}_{B})}
 
 \newcommand{\HBhN}{H^1_c(\Gamma_{0}(Np),\mathcal{D}_{B_{h,N}})}
 \newcommand{\HL}{H^1_c(\Gamma_{0}(Np^{m}),L_{k,\chi}(R_{\kappa}))}
 \newcommand{\HLs}{H^1_c(\Gamma_{0}(Np^{m}),L_{2k,\chi^2}(R_{\kappa}))}
 
 \newcommand{\HLsp}{H^1_{par}(\Gamma_{0}(Np^{m}),L_{2k,\chi^2}(\bC_p))}
 \newcommand{\HR}{H^1_c(\Gamma_{0}(Np),\mathcal{D}_{\mathcal{R}_{h,N}})}
 \newcommand{\ONc}{\Omega_{h,N}^{classical}}
 \newcommand{\RN}{\mathcal{R}_{h,N}}
 \newcommand{\RNT}{\tilde{\mathcal{R}}_{h,N}}
 \newcommand{\BN}{B_{h,N}}

 \newcommand{\Cusp}{S_{2k+2}(\Gamma_0(Np^m),\chi^2)}
 \newcommand{\Hcusp}{S_{k+\frac{3}{2}}(\Gamma_0(4Np^m),\chi^{'})}

 \newcommand{\DN}{\tilde{\mathcal{D}}(\mathbb{Z}_{p,N}^{\times})}
 
 \newcommand{\ZZ}{\mathbb{Z}_{p}^{\times} \times \mathbb{Z}_{p}}
 \newcommand{\TAR}{\otimes_{A(B_{h,N})} \mathcal{R}_{h,N}}
 \newcommand{\TNS}{\otimes_{_{\mathcal{D}(\mathbb{Z}_{p,N}^\times),\sigma}}}
 
 \newcommand{\NS}{\mathcal{D}(\mathbb{Z}_{p,N}^\times)}
 \newcommand{\NSN}{\mathcal{D}(\mathbb{Z}_{p}^\times)[\Delta_N]}
 \newcommand{\ten}{\otimes_{_{\mathcal{D}(\mathbb{Z}_{p}^\times)}}}
 
 \newcommand{\Gam}{\Gamma_{0}(Np)}

\begin{document}

\title[Overconvergent Shintani Lifting for Coleman's $p$-adic family]
 {$P$-adic family of half-integral weight modular forms via
overconvergent Shintani lifting}

\author{Jeehoon Park}

\address{Jeehoon Park: Department of Mathematics, Boston Univ., Boston, MA,
02215, USA} \email{pineon@math.bu.edu}


\thanks{}

\thanks{}



\date{}

\dedicatory{}

\commby{}


\begin{abstract}
The classical Shintani map (see \cite{Shn}) is the
Hecke-equivariant map from the space of cusp forms of integral
weight to the space of cusp forms of half-integral weight. In this
paper, we will construct a Hecke-equivariant overconvergent
Shintani lifting which interpolates the classical Shintani lifting
$p$-adically, following the idea of G. Stevens in \cite{St1}. In
consequence, we get a formal $q$-expansion $\Theta$ whose
$q$-coefficients are in an overconvergent distribution ring, which
can be thought of $p$-adic analytic family of overconvergent
modular forms of half-integral weight, since the specializations
of $\Theta$ at the arithmetic weights are the classical cusp forms
of half-integral weight.

\end{abstract}

\maketitle

\tableofcontents

\subsection{Introduction.}\

The space of modular forms of half-integral weight (see
\cite{Shm2} for the precise definition) have been studied by many
people. One of the main reasons for this would be the facts that
their fourier coefficients are rich in arithmetic. For example,
the representation numbers of positive definite integral quadratic
forms of odd dimensions arise as fourier coefficients of theta
functions of half-integral weights. H. Cohen showed that his
generalized class numbers are the fourier coefficients of
Eisenstein series of half-integral weight (see \cite{Coh}).
Moreover, J. Waldspurger had shown (see \cite{Wa}) that the square
of the fourier coefficients of new cusp forms of half integral
weight (see \cite{Koh} for the theory of new forms in the
half-integral case) are closely connected to the special values of
$L$-function of corresponding integral weight cusp forms under the
Shimura correspondence.\\

It was Shimura who initiated the program which relates modular
forms of integral weight with modular forms of half-integral
weight. He defined the Hecke actions on half-integral weight
modular forms and constructed integral weight modular forms from
half-integral modular forms by using Hecke operators and
functional equations of Dirichlet L-functions associated to
concerning half-integral and integral weight modular forms, which
comes from explicit transformation formulas of some theta
functions. After Shimura's work, Shintani (\cite{Shn}) gave an
inverse construction to the Shimura's map, where he used the Weil
representation to construct the half-integral weight cusp form
$\theta(f)$ from an integral weight cusp form $f$. Furthermore, he
expressed the fourier coefficients of $\theta(f)$ as some period
integrals of $f$, hence can be described in terms of a certain
cohomology class associated to $f$ as in Shimura's work (see the
chapter 8 of \cite{Shm1}). G. Stevens realized that if we have a
$\Lambda$-adic cohomology class which interpolates the
p-stabilized ordinary newforms then it gives rise to a
$\Lambda$-adic version of Shintani lifting for Hida's universal
ordinary $\Lambda$-adic modular forms (see \cite{St1}). The main
goal of this paper is to generalize his result to the finite slope
non-ordinary case, i.e. Robert Coleman's $p$-adic analytic family
of overconvergent cuspidal eigenforms of finite slope.\\

In section 1, we recall the definition of overconvergent modular
symbols which is a compact supported group cohomology of a
congruence subgroup of $SL_2(\bZ$) with values in overconvergent
distribution and explain how it interpolates the classical cusp
forms by the specialization maps at arithmetic weights, i.e. how a
cohomological approach to $p$-adic family of modular forms works.
The more precise dictionary between Robert Coleman's $p$-adic
family of overconvergent eigenforms and overconvergent Hecke
eigensymbol (see the theorem \ref{Hecke eigensymbol exists}) will
be given in section 3. Note that before the final section 3, we
won't concentrate on Hecke eigenforms; instead, we deal with
modular forms which are not necessarily Hecke eigenforms.\\

In section 2, we give a purely cohomological formulation of the
Shintani lifting following G. Stevens' work (\cite{St1}), after
reviewing the Shimura correspondence, especially the Shintani's
construction of half-integral weight cusp forms from integral
weight cusp forms. Then $p$-adic lifting of the cohomological
Shintani lifting will be given (see the definition \ref{maindef}),
which we call overconvergent Shintani lifting. Its $p$-adic
interpolation property essentially asserts that the following
diagram commutes:

\begin{displaymath}
\xymatrix{ \Hc \ar[r]^{\ \ \ \ \Theta} \ar[d]^{\phi_{\kappa,*}} &
\DN[[q]] \ar[d]^{T_p^{m-1}\circ\tilde{\kappa}}\\
\HLs \ar[r]^{\ \ \ \ \ \ \ \ \ \ \  \Theta_{k,\chi}} &
\cR_{\kappa}[[q]] }
\end{displaymath}

where all the notations will be provided in the paper.
Furthermore, we describe how the natural Hecke action on
overconvergent modular symbol is transformed under the
overconvergent Shintani lifting (see the theorem \ref{main1}). The
main theorem of section 2 \ref{main1} implies that we can define
Hecke actions (see the definitions (\ref{Hecke action1}) and
(\ref{Hecke action2})) on the space of formal $q$-expansions with
coefficients in an overconvergent distribution ring so that the
overconvergent Shintani lifting is Hecke-equivariant. The
important feature of the cohomological (overconvergent) Shintani
lifting is that it's purely cohomological and algebraic which
depends only on the
arithmetic of integral indefinite binary quadratic forms.\\

In section 3, we state and prove our main theorem (the theorem
\ref{thm8}) which asserts the existence of formal $q$-expansion
$\Theta \in \RNT[[q]]$ whose specializations to classical
arithmetic points are, up to multiplication by scalars, the
$q$-expansions of Hecke eigenforms of half-integral weight (in the
sense of Shimura \cite{Shm2}), which is obtained via the classical
Shintani's map from the cuspidal Hecke eigenforms of tame
conductor $N$ with finite slope $\leq h$. In contrast to the
ordinary case, the existence of Hecke eigensymbol in $\Hc$ over
the whole weight space $\cX_N$ is not guaranteed in the
non-ordinary situation. But if we shrink the weight space small
enough, then we can prove the existence of Hecke eigensymbol over
that small affinoid domain. More precisely, G. Stevens proved that
there exists an affinoid subdomain $\BN$ in $\cX_N$ such that
$\HBhN$ admits a slope $\leq h$ decomposition (see the theorem
\ref{thmlocal}) and we prove the existence of Hecke eigensymbol in
the slope $\leq h$ part $\HBhN^{(\leq h)}$ (see the theorem
\ref{Hecke eigensymbol exists}). After introducing overconvergent
$p$-adic Hecke algebra for the slope $\leq h$ part $\HBhN^{(\leq
h)}$, we define a local overconvergent Shintani lifting over $\BN$
and prove its Hecke-equivariance and the main theorem as a simple
application of tools developed in the paper.\\

It should be mentioned that H. Hida constructed a $\Lambda$-adic
Shimura lifting in \cite{Hi2} which is an inverse to G. Stevens'
$\Lambda$-adic Shintani lifting (in \cite{St1}) which we are
generalizing to the finite slope non-ordinary case in our paper.
Also, Nick Ramsey constructed half-integral weight eigencurve and
studied overconvergent Shimura lifting in \cite{Ram1} and
\cite{Ram2}. The overconvergent Shintani lifting here can be
thought of an inverse to overconvergent Shimura lifting. It would
be worthwhile to write down a local (possibly global) rigid
analytic map from Coleman-Mazur integral weight eigencurve (see
\cite{Col-Mz} for the tame level 1 and \cite{Buz} for the higher
tame level $N$) to half-integral weight eigencurve (see
\cite{Ram2}), in the rigid analytic geometric language used in N.
Ramsey's paper \cite{Ram2}, though, in some sense, our overconvergent
Shintani lifting answers what it should be.\\

It's worthwhile to note that the recent work of Bertolini, Darmon,
and Tornaria about Hida families and Shimura lifts (see
\cite{BDT}) tells us that the $p$-adic derivative of a certain
linear combination of the fourier coefficients of a $\Lambda$-adic
cusp form of half-integral weight, evaluated at some arithmetic
point, are closely connected to the Stark-Heegner points on the
elliptic curve attached to a weight 2 cusp form which moves
$p$-adically in the corresponding $\Lambda$-adic cusp form of
integral weight under $\Lambda$-adic Shimura lifting. The fourier
coefficients of the universal overconvergent half-integral weight
modular forms (see the definitions \ref{univmod} and
\ref{univmodB}) might subject to the similar story so that the
$p$-adic derivative of their certain linear combination could be
closely related to the Heegner cycles in the non-ordinary and
higher weight $> 2$ case. It would certainly be a good project to
develop an analogous
theory in the non-ordinary and higher weight $>2$ case, following them.\\

The author thanks his thesis advisor G. Stevens for suggesting
this problem and for many invaluable conversations and help on the
subject, and Robert Pollack for helpful comments.

\section{Cohomological approach to Coleman's $p$-adic family of overconvergent modular forms}

Throughout the paper, we fix a prime $p \geq 5$ and a positive
integer $N$ with $(p,N)=1$; the "quadratic residue symbol"
$(\frac{a}{b})$ has the same meaning as in \cite{Shm2}. We use the
notations $\Delta_N$ and $(\bZ/N\bZ)^{\times}$ interchangeably.
Let $K$ be any finite extension field of $\bQ_p$ and $|\cdot|$ be
the complete non-archimedean absolute value of $K$ normalized by
$|p|=p^{-1}$. Let $\bC_p$ denote the $p$-adic completion of an
algebraic closure of $\bQ_p$. We fix, once and for all, a field
isomorphism between $\bC$ and $\bC_p$. Robert Coleman constructed
a $p$-adic family of overconvergent modular forms in \cite{Col1}.
It seems to be hard at the moment to define an overconvergent
Shintani lifting directly from the Coleman's way of describing the
$p$-adic family. But if we approach Coleman's $p$-adic family in a
cohomological way (using overconvergent modular symbols) as
developed by G. Stevens, then there is a natural way to get an
overconvergent Shintani lifting which interpolates the classical
Shintani $\theta$-lifts in \cite{Shn}. We start from a locally
analytic distribution space which will be the coefficient ring of
our compact supported cohomology for an arithmetic group and then
we describe overconvergent modular symbols with Hecke action on
them and how they interpolate classical modular forms by the
specialization maps.


\subsection{Overconvergent modular symbols}\

The goal of this section is to introduce overconvergent modular
symbols and the Hecke action on them. Let $S$ be a compact subset
of $\bQ_p$ or $\bQ_p \times \bQ_p$. We define $\cA(S)$ as the
space of $K$-valued locally analytic functions on $S$. We also
define $\cD(S)$ as the space of $K$-valued locally analytic
distributions on $S$, i.e. $\cD(S):= \textrm{Hom}_{cts}(\cA(S),K)$
with the strong topology. We will use the notation $\cA$ and
$\cD$, when $S=\ZZ$, i.e. $\cA=\cA(\ZZ)$ and $\cD=\cD(\ZZ)$. Then
$\cA$ is a reflexive and complete Hausdorff locally convex
$K$-vector space whose strong dual $\cD$ is a Fr\'{e}chet space.
We refer to
$\mathcal{x} 16$ of \cite{Sch} for the detailed definitions. \\

For any positive integer $M$, let
\begin{eqnarray*}
\Gamma_{0}(M):=\{{\left(%
\begin{array}{cc}
  a & b \\
  c & d \\
\end{array}%
\right) \in SL_2(\bZ) \mid c \equiv 0 \ (\textrm{mod} M)}\}
\end{eqnarray*}
 be the congruence subgroup of $SL_2(\bZ)$. Then
$\Gamma_0(Np)$ acts on $\ZZ$ by matrix multiplication on the
right, viewing the elements of $\ZZ$ as row vectors. So it induces
a left action on $\cA$ by $(\gamma \cdot f)(x,y):=f((x,y) \cdot
\gamma)$ where $(x,y)\in \ZZ$ and $f\in \cA$. Consequently it
induces a right action of $\Gamma_0(Np)$ on $\cD$ which is
uniquely characterized by the following integration formula:
\begin{eqnarray}
\int_{\ZZ}f(x,y) d(\mu|\gamma)(x,y)=\int_{\ZZ} (\gamma \cdot
f)(x,y) d\mu(x,y).\label{action1}
\end{eqnarray}
We refer to this action on $\cD$ as the dual action induced from
action on $\cA$. The scalar action of $\bZ_p^{\times}$ on $\ZZ$
induces a left action of $\cD(\bZ_p^{\times})$ on $\cA$ by
\begin{eqnarray}
(\nu \cdot f)(x,y):=\int_{\bZ_p^{\times}} f(x\cdot \lambda,y\cdot
\lambda) d\nu(\lambda)
\end{eqnarray}
for $\nu \in \cD(\bZ_p^{\times})$ and $f\in \cA$.

Hence this induces a right action of $\cD(\bZ_p^{\times})$ on
$\cD$ by the dual action, i.e. by the following formula:
\begin{eqnarray}
\int_{\ZZ}f(x,y) d(\mu|\nu)(x,y)=\int_{\ZZ} (\nu \cdot f)(x,y)
d\mu(x,y)
\end{eqnarray}
for $\mu \in \cD$ and $\nu \in \cD(\bZ_p^{\times})$. The action of
$\Gamma_0(Np)$ clearly commutes with the action of
$\cD(\bZ_p^{\times})$. Henceforth we may consider $\cD$ as a
$\cD(\bZ_p^{\times})[\Gamma_0(Np)]$-module.\\

Now we define the group of $\cD$-valued modular symbols over
$\Gamma_0(Np)$ whose elements we will call overconvergent modular
symbols (over $\Gamma_0(Np)$).
\begin{defn}
We define the group of $\cD$-valued modular symbol over
$\Gamma_0(Np)$ to be
\begin{eqnarray}
\textrm{Symb}_{\Gamma_0(Np)}(\cD):=\textrm{Hom}_{\Gamma_0(Np)}(\Delta_0,\cD)
\end{eqnarray}
where $\Delta_0:=Div^0(\mathbb{P}^1(\bQ))$ is the group of
divisors of degree 0 supported on the rational cusps
$\mathbb{P}^1(\bQ)$ of the upper half plane. \label{def1}
\end{defn}

In other words, $\Phi \in \textrm{Symb}_{\Gamma_0(Np)}(\cD)$ is an
additive homomorphism $\Phi:\Delta_0 \to \cD$ for which
$\Phi|\gamma=\Phi$ for all $\gamma \in \Gamma_0(Np)$ where the
action of $\gamma \in \Gamma_0(Np)$ on
$\textrm{Hom}(\Delta_0,\cD)$ is given by
\begin{eqnarray}
(\Phi|\gamma)(D):=\Phi(\gamma \cdot D)|\gamma \label{action}
\end{eqnarray}
for $D \in \Delta_0$.\\

We refer to \cite{Shm1} or the appendix of \cite{Hi} for the
definition of compactly supported
cohomology $H^1_c(\Gamma_0(Np),\cD)$.\\

It turns out there exists a canonical isomorphism (see \cite{A-S1}
for the proof)
\begin{eqnarray}
H^1_c(\Gamma_0(Np),\cD) = \textrm{Symb}_{\Gamma_0(Np)}(\cD).
\end{eqnarray}
Henceforth we will use the terms $\cD$-valued degree one compactly
supported cohomology and overconvergent modular symbol
interchangeably.\\

Consider the natural map
\begin{eqnarray*}
H^1_c(\Gamma_0(Np),\cD) \to H^1(\Gamma_0(Np),\cD)
\end{eqnarray*}
defined by sending a overconvergent modular symbol $\varphi$ to
the cohomology class represented by the 1-cocycle $\gamma \mapsto
\varphi(\{\gamma x\}-\{x\})$ for any fixed $x \in
\mathbb{P}^1(\bQ_p)$ (it is independent of the choice of fixed
element $x$). Then the parabolic cohomology (see [Shm1] or the
appendix of [Hi] for the definition of parabolic cohomology) is
canonically isomorphic to the image of the above natural map:
\begin{eqnarray}
H^1_c(\Gamma_0(Np),\cD) \to
H_{par}^1(\Gamma_0(Np),\cD)\hookrightarrow H^1(\Gamma_0(Np),\cD).
\label{sespar}
\end{eqnarray}

\begin{rmk}
If we replace $Np$ and $\cD$ by any positive integer $M$ and any
right $\bZ[\Gamma_0(M)]$-module respectively, (1.4)-(1.7) still
remains valid. \label{rem1}
\end{rmk}

We now give the Hecke-module structure of these cohomology groups.
Let $S_0(M)$ denote the semi-group

\begin{eqnarray*}
S_{0}(M):=\{{\left(%
\begin{array}{cc}
  a & b \\
  c & d \\
\end{array}%
\right) \in M_2^+(\bZ) \mid  c \equiv 0 \ (\textrm{mod} M)},\
(a,M)=1\}
\end{eqnarray*}
where $M$ is a positive integer and $M_2^+(\bZ)$ is the semi-group
of integral $2 \times 2$ matrices with positive determinant. Let
$R(\Gamma_0(M), S_{0}(M))$ be the $\bZ$-module generated by the
double cosets $\Gamma_0(M) \alpha \Gamma_0(M), \alpha \in S_0(M)$.
This can be equipped with a ring structure by defining
multiplication between two double cosets. $R(\Gamma_0(M),
S_{0}(M))$ with this multiplication law extended $\bZ$-linearly,
becomes a commutative ring with $\Gamma_0(M)=\Gamma_0(M) \cdot 1
\cdot \Gamma_0(M)$ as the unit element. It is called a Hecke ring
with respect to a
Hecke pair $(\Gamma_0(M), S_0(M))$.\\

For each positive integer $n$, denote by $T_n$ the formal sum of
all double cosets $\Gamma_0(M) \alpha \Gamma_0(M)$ with $\alpha
\in S_0(M)^{n}:=\{\alpha \in S_0(M) \mid det(\alpha)=n\}$ in
$R(\Gamma_0(M), S_{0}(M))$. For example, $T_p=\Gamma_{0}(M)\left(%
\begin{array}{cc}
  1 & 0 \\
  0 & p \\
\end{array}%
\right) \Gamma_{0}(M)$ for every prime $p$. For primes $p$
satisfying $(p,M)=1$, let
\begin{eqnarray}
T_{p,p}=\Gamma_0(M)\left(
                     \begin{array}{cc}
                       p & 0 \\
                       0 & p \\
                     \end{array}
                   \right)
\Gamma_0(M).\label{Tpp}
\end{eqnarray}

The structure of $R(\Gamma_0(M), S_{0}(M))$ is given by the
following theorem [Shm1, Theorem 3.34]:

\begin{prop}
(1) $R(\Gamma_0(M), S_{0}(M))$ is a polynomial ring over $\bZ$ in
the variables $T_{p,p}$ for all primes $p$ not dividing $M$ and
$T_p$ for all primes $p$.

(2) $R(\Gamma_0(M), S_{0}(M))\otimes_{\bZ} \bQ_p$ is a polynomial
ring over $\bQ_p$ in the variables $T_n$ for all $n>0$, i.e.
\begin{equation}
R(\Gamma_0(M), S_{0}(M))\otimes_{\bZ} \bQ_p=\bQ_p[T_n : n \in \bN]
\end{equation} \label{prop1}
\end{prop}

We will call $R_{\bQ_p}:=R(\Gamma_0(M), S_{0}(M))\otimes_{\bZ}
\bQ_p$ the
abstract Hecke algebra over $\bQ_p$ for $(\Gamma_0(M),S_0(M))$. \\
If $p \mid M$, then we use the notation $U_p$ instead of $T_p$. If
$A$ is a right $S_0(M)$-module, then $A^{\Gamma_0(M)}$ has a
natural action by the Hecke operators $T_n$ induced by, for $a \in
A^{\Gamma_0(M)}$ and $\alpha \in S_0(M)^n$,
\begin{equation}
a|\Gamma_{0}(M) \alpha \Gamma_{0}(M):= \sum_i a|\alpha_i,
\end{equation}
where $ \{\alpha_i\}$ is the set of representatives for
$\Gamma_{0}(M) \backslash \Gamma_{0}(M) \alpha \Gamma_{0}(M)$.
Since the action of $\Gamma_0(M)$ on $\cD$ extends to an action of
$S_0(M)$ on $\cD$ by the same formula (\ref{action1}),
$\textrm{Hom}(\Delta_0,\cD)$ is a right $S_0(M)$-module. Therefore
$H_c^1(\Gamma_0(M),\cD)=\textrm{Hom}_{\Gamma_0(M)}(\Delta_0,\cD)$
can be viewed as $R_{K}$-module by the action described above,
where $R_{K}=R(\Gamma_0(M), S_{0}(M))\otimes_{\bZ}{K}$. In
particular, if $p \mid M$, then $U_p$ acts on
$H^1_c(\Gamma_{0}(M),\cD)$ by
\begin{equation}
\Phi|U_p = \sum_{i=1}^{p} \Phi|\left(%
\begin{array}{cc}
  1 & i \\
  0 & p \\
\end{array}%
\right),
\end{equation}
for $\Phi \in H^1_c(\Gamma_{0}(M),\cD)$.\\

Because $\cD$ has an action of $S_0(M)$ by the same formula
(\ref{action1}), we can also define the $R_{K}$-actions on
$H^1(\Gamma_0(M),\cD)$ and $H_{par}^1(\Gamma_0(M),\cD)$. In fact
the maps in the sequence (\ref{sespar}) are $R_{K}$-module
homomorphisms, i.e. Hecke-equivariant homomorphism.


The matrix $\iota:=\left(%
\begin{array}{cc}
  1 & 0 \\
  0 & -1 \\
\end{array}%
\right)$ induces natural involutions on $H^1_c(\Gamma_0(M),\cD)$
and $H^1(\Gamma_0(M),\cD)$. The involution $\iota$ decomposes each
of cohomology groups as $\pm$ eigenspaces:
\begin{eqnarray}
H=H^+ \oplus H^-.
\end{eqnarray}
We know that each cohomology class $\varphi$ decomposes as
$\varphi=\varphi^+ +\varphi^-$ where
$\varphi^{\pm}:=\frac{1}{2}(\varphi {\pm} \varphi|\iota)$ and
$\varphi^{\pm}|\iota={\pm}\varphi$.\\

Let $R$ be any commutative ring with unity. We recall the
definition of $L_{k,\chi}(R)$ and $L_{k,\chi}^{\ast}(R)$, where
$\chi$ is an $R$-valued Dirichlet character of
$(\bZ/M\bZ)^{\times}$ for some positive integer $M$, from 4.1 of
\cite{St1}. For any integer $k\geq 0$ we define $Sym^k(R^2)$
(respectively, $Sym^k(R^2)^{\ast}$) as the free $R$-module
generated by the divided powers monomials $\frac{X^{n}}{n!}\cdot
\frac{Y^{k-n}}{(k-n)!}$(respectively, the monomials $X^n \cdot
Y^{k-n}$) for $0\leq k \leq n$. Then $L_{k,\chi}(R)$
(respectively, $L_{k,\chi}^{\ast}(R)$) is $R[\Gamma_0(M)]$-module
whose underlying $R$-module is $Sym^k(R^2)$ (respectively,
$Sym^k(R^2)^{\ast}$) equipped with the following
action of $S_0(M)$: for $\gamma=\left(%
\begin{array}{cc}
  a & b \\
  c & d \\
\end{array}%
\right)\in S_0(M)$,
\begin{displaymath}
(F|\gamma)(X,Y):= \left\{ \begin{array}{ll}
\chi(a)\cdot F((X,Y)\cdot\gamma^{*}) & \textrm{for $F \in L_{k,\chi}(R)$}\\
\\
\chi(d)\cdot F((X,Y)\cdot\gamma^{*}) & \textrm{for $F \in
L_{k,\chi}(R)^{\ast}$}
\end{array} \right.
\end{displaymath}where $\gamma^{*}=\left(
                    \begin{array}{cc}
                      d & -b \\
                     -c &  a \\
                    \end{array}
                  \right)$.
 There is a unique perfect $R$-linear pairing
$\langle\cdot,\cdot\rangle:Sym^k(R^2) \times Sym^k(R^2)^{\ast} \to
R$ satisfying $\big\langle \frac{X^iY^{k-i}}{i!(k-i)!},X^{k-j}Y^j
\big\rangle=(-1)^j\delta_{ij}$ where $\delta_{ij}$ is the
Kronecker delta. In particular, we have
\begin{eqnarray}
\big\langle \frac{(aY-bX)^k}{k!},P(X,Y)
\big\rangle=P(a,b)\label{pairing}
\end{eqnarray}
for every $(a,b)\in R^2$ and $P\in Sym^k(R^2)^{\ast}$.\\

Let $S_{k+2}(\Gamma_0(M),\chi)$ (respectively
$M_{k+2}(\Gamma_0(M),\chi)$) be the space of weight $k+2$
cusp(respectively modular) forms of level $\Gamma_0(M)$ with
nebentype character $\chi$ and $S_{k+2}(\Gamma_0(M),\chi;R)$ the
subspace of $S_{k+2}(\Gamma_0(M),\chi)$ whose Fourier coefficient
lie in $R$ where $R$ is a subfield of $\bC$. We recall the
Eichler-Shimura isomorphism (see [Shm1])
\begin{eqnarray*}
S_{k+2}(\Gamma_0(M),\chi) &{}^{\sim \atop \To}&
H^1_{par}(\Gamma_0(M),L_{k,\chi}(\bC))^{\pm}.\\
f \ \ \  &\mapsto& \ \ \ \Psi_f^{\pm}
\end{eqnarray*}
for a positive integer $M$ and Dirichlet character $\chi$ defined
modulo $M$.

By the Manin-Drinfeld principle, there is a Hecke-equivariant
section to the natural map (\ref{sespar}):
\begin{eqnarray*}
H^1_{par}(\Gamma_0(M),L_{k,\chi}(R)) &\to&
H^1_c(\Gamma_0(M),L_{k,\chi}(R)), \\
\Psi \ \ \ &\mapsto& \ \ \ \psi
\end{eqnarray*}
for a positive integer $M$ and Dirichlet character $\chi$ defined
modulo $M$, if $R$ is the field of characteristic 0. Henceforth we
can consider the element of parabolic cohomology
$H^1_{par}(\Gamma_0(M),L_{k,\chi}(R))$ as a modular symbol with
values in $L_{k,\chi}(R)$.

\subsection{$p$-adic interpolation of overconvergent modular
symbols}\

In this section we review how overconvergent modular symbols can
be understood as $p$-adic families of modular forms (see
\cite{St3}). For that we recall the definition of (arithmetic)
weights and introduce the specialization map for an arithmetic weight.\\

We refer to an element of
$\textrm{Hom}_{cts}(\bZ_{p,N}^{\times},\bC_p^{\times})$ as a
weight. For $t \in \bZ_{p,N}^{\times}$ we let $t_p$ and $t_N$ be
the coordinates of $t$ under the canonical decomposition
$\bZ_{p,N}^{\times}=\bZ_{p}^{\times}\times(\bZ/N\bZ)^{\times}$.
Let $\Lambda_N=\bZ_p[[\bZ_{p,N}^{\times}]]$. Then there is a
natural identification between
\begin{eqnarray}
\textrm{Hom}_{cts}(\Lambda_N,\bC_p) &\textrm{and}&
\textrm{Hom}_{cts}(\bZ_{p,N}^{\times},\bC_p^{\times}).
\label{points}
\end{eqnarray}
where the first Hom denotes continuous algebra homomorphisms and
second denotes continuous group homomorphisms.

\begin{defn}
A character(weight) $\kappa \in
\textrm{Hom}_{cts}(\bZ_{p,N}^{\times},\bC_p^{\times})$ is called
arithmetic of signature of $(k,\chi)$, if it satisfies
\begin{eqnarray*}
\kappa(t)=\chi(t)\cdot t_p^k
\end{eqnarray*}
for some $\chi$, a finite order character of $\bZ_{p,N}^{\times}$,
and $k\in \bZ^{\geq0}$. \label{def3}
\end{defn}

Note that we have a natural isomorphism between
$\cD(\bZ_{p,N}^{\times})$ and $\NSN$. Accordingly, we refer to
them interchangeably. Let's define
\begin{eqnarray*}
\cD_N:=\cD\ten \NS,
\end{eqnarray*}
where $\gamma \in S_0(Np)$ acts on $\cD_N$ as follows:
\begin{eqnarray}
 (\mu \ten \nu)\mid(\gamma):= \mu \mid \gamma \ten
 [a(\gamma)]_N\cdot \nu \ \ , \label{actionN}
\end{eqnarray}
for $\mu \in \cD$, $\nu \in \NS$ and $[a(\gamma)]_N \in \Delta_N \
( \textrm{here}\ a(\gamma)\ \textrm{is the upper left
 entry of $\gamma$}$).
So $\cD_N$ can be thought as $\NS[\Gam]$-module and $\Hc$ is
defined and has a natural $\NS$-module structure.

\begin{defn}
We define the abstract overconvergent Hecke algebra of tame
conductor $N$ to be the free polynomial algebra
\begin{eqnarray*}
\cH:=\NS[T_n :  n \in \bN]
\end{eqnarray*}
generated by over $\NS$ by $T_n \in R(\Gam,S_0(Np))$ for $n\in
\bN$.\label{abstract Hecke algebra}
\end{defn}
Note that $\cH \cong R(\Gam,S_0(Np))\otimes_{\bZ}\NS$. If $\kappa
\in \textrm{Hom}_{cts}(\bZ_{p,N}^{\times},K^{\times})$ is any
arithmetic point of signature $(k, \chi)$, then let
$R_\kappa:=\kappa(\NS)$ where $\kappa(r):=r(\kappa)$ for
$r\in\NS$. We can define a $R_\kappa$-linear map, called a
specialization map for $\kappa$,
\begin{eqnarray}
\phi_\kappa: \cD_N=\cD\ten \NS \to L_{k,\chi}(R_\kappa) \label{sp}
\end{eqnarray}
by
\begin{eqnarray}
\phi_\kappa(\mu \otimes
r):=\kappa(r)\cdot\int_{\ZZ}\chi_p(x)\frac{(xY-yX)^k}{k!}d\mu(x,y)
\label{specialize}
\end{eqnarray}
for $\mu \in \cD$ and $r \in \NS$, where we factor $\chi=\chi_N
\cdot \chi_p$ with $\chi_N$ defined modulo $N$ and $\chi_p$
defined modulo a power of $p$ and $\kappa(r):=r(\kappa)$. A simple
computation confirms that if $\chi$ is defined modulo $Np^m$ then
$\phi_\kappa$ is a $\Gamma_0(Np^m)$-module homomorphism, and hence
induces a homomorphism on cohomology groups:
\begin{eqnarray}
\phi_{\kappa,\ast}:\Hc \to
H^1_c(\Gamma_0(Np^m),L_{k,\chi}(R_\kappa)).\label{speci}
\end{eqnarray}
The map $\phi_{\kappa,\ast}$ is Hecke-equivariant. Since
1-cocycles in $\HL$ can be viewed as classical modular forms via
the Manin-Drinfeld principle and the Eichler-Shimura isomorphism,
we can interpret an overconvergent modular symbol as a $p$-adic
analytic family of modular forms by varying $k$ and $\chi$.

\section{Overconvergent Shintani Lifting}
We first review the original Shintani lifting in \cite{Shn} and
its cohomological realization in \cite{St1}. Then we give an
overconvergent version of these constructions. We will prove the
Hecke-equivariance of our overconvergent Shintani lifting and
decribe the Hecke action explicitly on formal $q$-expansions with
coefficient in $\DN:=\NS\TNS \NS$ (see (\ref{metaplectic})), which
is given as the image of the overconvergent Shintani lifting. This
formal $q$-expansion can be viewed as the $p$-adic family of
overconvergent half-integral weight modular forms.

\subsection{Classical Shintani lifting and its cohomological
interpretation}\

We briefly summarize Shintani's construction in \cite{Shn} of a
Hecke-equivariant map from integral weight cusp forms to
half-integral weight cusp forms. We then review its cohomological
version which is well-suited to $p$-adic variation of Shintani
map, following \cite{St1} whose results we generalize to the
non-ordinary case. Let $Q=Q(X,Y)=aX^2+bXY+cY^2$ be an integral
(meaning $a,b,c \in \bZ$) binary quadratic form. We define the
discriminant of $Q$ as $\delta_Q:=b^2-4ac$. We call $Q$ an
indefinite quadratic form if $\delta_Q > 0.$ Let $\cF$ be the
space of integral indefinite binary quadratic forms.

\begin{defn}
For $M \in \bN$, we define the set $\mathcal{F}_M$ as follows:
\begin{eqnarray*}
\cF_{M}:=\left\{ \begin{array}{ll} \{Q(X,Y)=aX^2+bXY+cY^2\in \cF
\mid \ (a,M)=1, \ M \mid b, M \mid
c\} &  \textrm{if $M$ is odd}\\
\\
\{Q(X,Y)=aX^2+bXY+cY^2\in \cF \mid\
(a,M)=1,\ 2M \mid b, M \mid c\} & \textrm{if $M$ is even}  \\
\end{array} \right.\label{def11}
\end{eqnarray*}

\end{defn}

If $Q(X,Y)=aX^2+bXY+cY^2 \in \cF$ satisfies $(a,b,c)=1$, we call
$Q$ a primitive quadratic form. The congruence subgroup
$\Gamma_0(M)$ acts on $\mathcal{F}_M$ by the following formula:
\begin{eqnarray}
(Q|\gamma)(X,Y):=Q((X,Y)\cdot\gamma^{-1})\label{actionQ}
\end{eqnarray}
for $Q\in \mathcal{F}_M$ and $\gamma\in \Gamma_0(M)$. It is easy
to check that this action of $\Gamma_0(M)$ preserves $\cF_M$. If a
binary quadratic form $Q$ is allowed to have rational
coefficients, then we define the action of $GL_2(\bQ)$ on $Q$ by
the same formula (\ref{actionQ}).
\\

Now we associate a pair of points $\omega_Q,\omega_Q'\in
\bP^1(\bR)=\bR \cup (i\infty)$ to each integral indefinite binary
quadratic form $Q \in \mathcal{F}_M$ following [Shn]:
\begin{displaymath}
(\omega_Q,\omega_Q') := \left\{ \begin{array}{ll}
\big(\frac{b+\sqrt{\delta_Q}}{2c},\frac{b-\sqrt{\delta_Q}}{2c}) & \textrm{if $c \neq 0$}\\
 & \\

 (i\infty,\frac{a}{b})& \textrm{if $c=0$ and $b>0$}\\
 & \\

(\frac{a}{b}, i\infty)& \textrm{if $c=0$ and $b<0$}
\end{array} \right.
\end{displaymath}

\begin{defn}
We define the oriented geodesic path $C_Q$ in the upper half plane
$\mathfrak{h}$ following \cite{Shn}:
\begin{displaymath}
C_Q := \left\{ \begin{array}{ll}
\textrm{the oriented geodesic path joining $\omega_Q$ to $\omega_Q'$} & \textrm{if $\delta_Q$ is a perfect square}\\
& \\
\textrm{the oriented geodesic path joining $\omega$ to
$\gamma_Q(\omega)$} & \textrm{otherwise},
\end{array} \right.
\end{displaymath}
\label{cycle}where $\omega$ is an arbitrary point in $\bP^1(\bQ)$ and $\gamma_Q=\left(%
\begin{array}{cc}
  r & s \\
  t & u \\
\end{array}%
\right)$ is the unique generator (satisfying $r-t\omega_Q >1$) of
$\Gamma_0(M)_Q^+$ which is the index 2 subgroup (consisting of
matrices with positive traces) of the stabilizer of $Q \in
\mathcal{F}_M$ in $\Gamma_0(M)$. See [Shn] and \cite{St1} for more
details.
\end{defn}

For a given Dirichlet character $\chi$ defined modulo $M$, we
define a new Dirichlet character $\chi'$ modulo $4M$ by
\begin{eqnarray}
\chi'(d):=\chi(d)\cdot\Big(\frac{(-1)^{k+1}M}{d}\Big), \ d
\in(\bZ/4M\bZ)^{\times}.
\end{eqnarray}

For each quadratic form $Q=Q(X,Y)=aX^2+bXY+cY^2 \in \mathcal{F}_M$
and a Dirichlet character $\chi$ defined modulo $M$, we define
\begin{eqnarray}
\chi(Q):=\chi(a).
\end{eqnarray}

Now we're ready to construct the classical Shintani
$\theta$-lifting. We define

\begin{displaymath}
\theta_{k,\chi}(f)=\theta_{k,\chi}(f,z) := \left\{
\begin{array}{ll} \sum_{Q\in \mathcal{F}_M/\Gamma_0(M)}
I_{k,\chi}(f,Q)q^{\delta_Q/M} & \textrm{if $M$ is odd,}\\
\\
\sum_{Q\in \mathcal{F}_M/\Gamma_0(M)}
I_{k,\chi}(f,Q)q^{\delta_Q/4M} & \textrm{if $M$ is even.}
\end{array} \right.
\end{displaymath}


for $f \in S_{2k+2}(\Gamma_0(M),\chi^2)$ and $z \in \mathfrak{h}$
(the upper half plane), where $q:=e^{2\pi iz}$ and the coefficient
in the $q$-expansion $I_{k,\chi}(f,Q)$ is defined by
\begin{eqnarray}
I_{k,\chi}(f,Q):= \chi(Q)\cdot \int_{C_Q}f(\tau)Q(1,-\tau)^k
d\tau.\label{Shintai integral}
\end{eqnarray}
The integral $I_{k,\chi}(f,Q)$ converges for cusp forms $f\in
S_{2k+2}(\Gamma_0(M),\chi^2)$ and depends only on the
$\Gamma_0(M)$-orbit of $Q$ in $\mathcal{F}_M$.\\

Recall that $S_{k+\frac{3}{2}}(\Gamma_0(4M),\chi')$ is the space
of cusp forms of level $4M$, weight $k+\frac{3}{2}$, and nebentype
character $\chi'$ (a Dirichlet character modulo $4M$). We refer to
\cite{Shm2} for detailed definitions of half-integral weight
modular forms. Shintani proved the following theorem (Theorem 2 in
[Shn]).

\begin{thm}
Let $k \geq 0$ and $\chi$ be a Dirichlet character defined modulo
$M$. Then for each $f \in S_{2k+2}(\Gamma_0(M),\chi^2)$, the
series $\theta_{k,\chi}(f,z)$ is the $q$-expansion of a
half-integral weight cusp form in
$S_{k+\frac{3}{2}}(\Gamma_0(4M),\chi')$. Moreover, the map
\begin{eqnarray}
\theta_{k,\chi}: S_{2k+2}(\Gamma_0(M),\chi^2) \to
S_{k+\frac{3}{2}}(\Gamma_0(4M),\chi')
\end{eqnarray}
is a Hecke-equivariant $\bC$-linear map, i.e.
\begin{eqnarray}
\theta_{k,\chi}(f|T_l)=\theta_{k,\chi}(f)|T_{l^2}
\end{eqnarray}
for any positive odd prime number $l$.\label{Shintani's main thm}
\end{thm}

We recall the definition of the Hecke operator $T_{p}$ on
$S_{k+\frac{3}{2}}(\Gamma_0(4M),\chi')$ when $p|M$. It is given on
$q$-expansions by
\begin{eqnarray}
\big(\sum_{n=1}^{\infty}\beta_n q^n
\big)|T_p:=\sum_{n=1}^{\infty}\beta_{pn} q^n.\label{Hecke p}
\end{eqnarray}

Proposition 1.5 of [Shm2] tells us that $T_p$ with $p|M$ induces a
map
\begin{eqnarray*}
T_p:S_{k+\frac{3}{2}}(\Gamma_0(4M),\chi')\to
S_{k+\frac{3}{2}}\big(\Gamma_0(4M),\chi'\cdot(\frac{p}{\cdot})\big)
\end{eqnarray*}
So $T_p$ preserves the level but multiplies the nebentype
character by the quadratic character $(\frac{p}{\cdot})$. For a
rational prime $l$ which doesn't divide $4M$, the Hecke operator
$T_{l^2}$ on $S_{k+\frac{3}{2}}(\Gamma_0(4M),\chi')$ preserves the
space $S_{k+\frac{3}{2}}(\Gamma_0(4M),\chi')$. Theorem 1.7 of
\cite{Shm2} tells us that
\begin{eqnarray}
\big(\sum_{n=1}^{\infty}\beta_n q^n
\big)|T_{l^2}:=\sum_{n=1}^{\infty}\big(\beta_{l^2n}+\chi'(l)\big(\frac{-1}{l}\big)^{k+1}\big(\frac{n}{l}\big)l^k
\beta_n + \chi'(l^2)l^{2k+1}\beta_{\frac{n}{l^2}}\big)
q^n.\label{Heckep}
\end{eqnarray}

Now let's turn to cohomological Shintani liftings. For each $Q\in
\mathcal{F}_M$, let
\begin{displaymath}
D_Q:=\partial C_Q= \left\{ \begin{array}{ll}
\{\omega_Q\}-\{\omega'_Q\} & \textrm{if $\delta_Q$ is a perfect square}\\
& \\ \{\gamma_Q(\omega)\}-\{\omega\} & \textrm{otherwise}.
\end{array} \right.
\end{displaymath}
Note that $D_Q \in \Delta_{0}$.
\begin{defn}
Let $R$ be a commutative $\bZ[\frac{1}{6}]$-algebra. Let $k\geq 0$
be an integer and $\chi$ be an $R$-valued Dirichlet character
defined modulo $M$. We define
\begin{eqnarray}
J_{k,\chi}(\phi,Q):=\chi(Q) \cdot \big{\langle} \phi(D_Q),Q^k
\big{\rangle} \in R
\end{eqnarray}
for $\phi \in
H^1_c(\Gamma_0(M),L_{2k,\chi^2}(R))=\mathrm{Hom}_{\Gamma_0(M)}(\Delta_0,L_{2k,\chi^2}(R))$
and $Q\in \mathcal{F}_M$ (so $Q^k\in Sym^k(R)^{*}$), where
$\phi(D_Q)$ is the value of any 1-cocycle representing $\phi$.
Note that the definition of $J_{k,\chi}(\phi,Q)$ is independent of
the choice of a representative cocycle and depends only on the
$\Gamma_0(M)$-orbit of $Q$ in $\mathcal{F}_M$.
\end{defn}

\begin{defn}
We define the cohomological Shintani lifting
$\Theta_{k,\chi}:H^1_c(\Gamma_0(M),L_{2k,\chi^2}(R)) \to R[[q]]$
by

\begin{displaymath}
\Theta_{k,\chi}(\phi):= \left\{
\begin{array}{ll} \sum_{Q\in \mathcal{F}_M/\Gamma_0(M)}
J_{k,\chi}(\phi,Q)q^{\delta_Q/M} \in R[[q]] & \textrm{if $M$ is odd,}\\
\\
\sum_{Q\in \mathcal{F}_M/\Gamma_0(M)}
J_{k,\chi}(\phi,Q)q^{\delta_Q/4M} \in R[[q]] & \textrm{if $M$ is
even.}
\end{array} \right.\label{cs1}
\end{displaymath}for each cohomology class $\phi \in
H^1_c(\Gamma_0(M),L_{2k,\chi^2}(R))$.
\end{defn}

Note that $\Theta_{k,\chi}$ is an $R$-linear map. The following
result was proved in \cite{St1} (Proposition (4.3.3)).

\begin{prop}
Let $R$ be a commutative $\bZ[\frac{1}{6}]$-algebra. Let $k\geq 0$
be an integer and $\chi$ be an $R$-valued Dirichlet character
defined modulo $M$.\\
(a) For $\varphi \in H^1_c(\Gamma_0(M),L_{2k,\chi^2}(R))$ and $\iota=\left(%
\begin{array}{cc}
  1 & 0 \\
  0 & -1 \\
\end{array}%
\right)$
\begin{eqnarray}
\Theta_{k,\chi}(\varphi|\iota)=-\Theta_{k,\chi}(\varphi).
\end{eqnarray}
(b) For $f \in S_{2k+2}(\Gamma_0(M),\chi^2)$ and $\psi_f\in
H^1_c(\Gamma_0(M),L_{2k,\chi^2}(\bC))$ which is the image of $f$
under the composition of the Eichler-Shimura map and the
Manin-Drinfeld section,
\begin{eqnarray}
\Theta_{k,\chi}(\psi_f)=\Theta_{k,\chi}(\psi_f^-)=\theta_{k,\chi}(f)
\end{eqnarray}
where $\theta_{k,\chi}(f)$ is defined in (2.4) and $\psi_f^-:=\frac{1}{2}\cdot(\psi_f-\psi_f|\iota)$.\\
(c) If $R$ is the field of characteristic 0 (for example $R=K$)
and $\varphi\in H^1_c(\Gamma_0(M),L_{2k,\chi^2}(R))$ is the image
of $\psi \in H^1_{par}(\Gamma_0(M),L_{2k,\chi^2}(R))$ under the
Manin-Drinfeld section, then
\begin{eqnarray}
\Theta_{k,\chi}(\varphi) \in S_{k+\frac{3}{2}}(\Gamma_0(4M),\chi';
R)
\end{eqnarray}
where $\chi'$ is defined in (2.2).\label{prop4}
\end{prop}


\subsection{Overconvergent Shintani Lifting and a universal overconvergent half-integral weight modular forms}\

In this section we define the Hecke-equivariant overconvergent
Shintani lifting (see \cite{St1} for $\Lambda$-adic Shintani
liftings) and describe the Hecke-action explicitly on the image of
the lifting which we define to be the universal overconvergent
half-integral weight modular form. Let $\sigma:\bZ_{p,N}^{\times}
\to \bZ_{p,N}^{\times}$ be the ring homomorphism sending $t$ to
$t^2$. We use the same notation for the ring homomorphism
\begin{eqnarray}
\sigma:\NS \to \NS \label{double}
\end{eqnarray}
induced from $\sigma:\bZ_{p,N}^{\times} \to \bZ_{p,N}^{\times}$,
i.e. $\sigma(\mu)(f(t))=\mu(f(t^2))$ for $\mu \in \NS$ and $f\in
\cA(\bZ_{p,N}^{\times})$. We define the space of metaplectic
locally analytic distributions on $\bZ_{p,N}^{\times}$ by
\begin{eqnarray}
\DN:=\NS \TNS \NS\label{metaplectic}
\end{eqnarray}where the tensor product is taken with respect to $\sigma:\NS \to
\NS$.
We regard $\DN$ as a $\NS$-algebra by the structure homomorphism
\begin{eqnarray}
\NS \to \DN \label{absm}
\end{eqnarray}
given by $r \mapsto r\TNS1$.
\\

The following definition is the key step to construct
overconvergent Shintani liftings.
\begin{defn}
For each $Q\in \cF_{Np}$ we define a map $J_{Q}:\cD \to
\cD(\bZ_p^{\times})$ by the following integration formula:
\begin{eqnarray}
\int_{\bZ_p^{\times}}f(z)
d\big(J_{Q}(\mu)\big)(z)=\int_{\bZ_p^{\times} \times
\bZ_p}f(Q(x,y)) d\mu(x,y)\label{JQmap}
\end{eqnarray}
for a given $\mu \in \cD$ and a locally analytic function $f\in
\cA(\bZ_p^{\times})$. This is well defined, since $Q(X,Y) \in
\bZ_p^{\times}$ for $Q\in \cF_{Np}$.
\end{defn}

We note that $J_Q$ is not a $\cD(\bZ_p^{\times})$-linear map. In
order to make it $\cD(\bZ_p^{\times})$-linear, we twist the target
space by $\sigma$ (see (\ref{metaplectic})). Then we extend $J_Q$
to $\tilde{J}_Q:\cD_N \to \DN$ by
\begin{eqnarray*}
\tilde{J}_Q(\nu \ten r):=1 \TNS J_{Q}(\nu)\cdot\sigma(r)\cdot[Q]_N
\in \DN
\end{eqnarray*}
for $\nu \ten r \in \cD \ten \NS=\cD_N$. Here $[Q]_N:=[a]_N \in
\Delta_N$ for each quadratic form $Q(X,Y)=aX^2+bXY+cY^2 \in
\mathcal{F}_{Np}$.

\begin{prop}
Let $Q\in \cF_{Np}$. The map $\tilde{J}_{Q}:\cD_N \to \DN$ is the
unique $\NS$-module homomorphism such that for all arithmetic
$\tilde{\kappa} \in
\textrm{Hom}_{cts}(\bZ_{p,N}^{\times},\bC_p^{\times})$ of
signature $(k,\chi)$ and the associated arithmetic $\kappa \in
\textrm{Hom}_{cts}(\bZ_{p,N}^{\times},\bC_p^{\times})$ by
$\sigma:\bZ_{p,N}^{\times} \to \bZ_{p,N}^{\times}$, i.e.
$\kappa:=\tilde{\kappa} \circ \sigma$ (so $\kappa$ has a signature
$(2k,\chi^2)$), we have
\begin{eqnarray}
\tilde{\kappa}(\tilde{J}_Q(\mu))=\chi(Q)\cdot\langle\phi_{\kappa}(\mu),Q^k\rangle
\label{int}
\end{eqnarray}
for $\mu \in \cD_N$. Here $\phi_{\kappa}$ is the specialization
map in (\ref{specialize}) and $\langle\cdot, \cdot \rangle$ is the
pairing defined by (\ref{pairing}). \label{JQ}
\end{prop}

\begin{proof}
Because of the construction, $\NS$-linearity is clear. If the
interpolation property (\ref{int}) at arithmetic points is proven,
then the uniqueness follows. So we concentrate on proving
(\ref{int}). Let $\mu=\nu \ten r \in \cD_N=\cD\ten \NS $, and let
$\tilde{\kappa} \in
\textrm{Hom}_{cts}(\bZ_{p,N}^{\times},\bC_p^{\times})$ of
signature $(k,\chi)$ associated to an arithmetic $\kappa \in
\textrm{Hom}_{cts}(\bZ_{p,N}^{\times},\bC_p^{\times})$ by
$\sigma:\NS \to \NS$ (so $\kappa$ has a signature $(2k,\chi^2)$).
We then calculate

\begin{eqnarray*}
\tilde{\kappa}(\tilde{J}_Q(\mu)) &=& \kappa(r)\cdot\tilde{\kappa}\big([Q]_N\cdot J_Q(\nu)\big)  \\
 &=&\tilde{\kappa}([Q]_N)\cdot \kappa(r) \cdot \int_{\bZ_p^{\times}}\tilde{\kappa}(z) dJ_{Q}(\nu)(z)\\
 &=&\chi_{N}(Q)\cdot \kappa(r)\cdot
 \int_{\bZ_p^{\times}\times\bZ_p}\tilde{\kappa}(Q(x,y))d\nu(x,y) \\
 &=& \chi_{N}(Q)\cdot \kappa(r)\cdot \int_{\bZ_p^{\times}\times \bZ_p}
 \chi_p(Q(x,y))\cdot Q(x,y)^k d\nu(x,y) \\
 &&\textrm{(because $Q(x,y) \in \bZ_p$ for $(x,y) \in \ZZ$)}\\
 &=& \chi_{N}(Q)\chi_{p}(a) \kappa(r) \cdot \int_{\bZ_p^{\times}\times\bZ_p}
 \chi_p(x^2)\Big\langle\frac{(xY-yX)^{2k}}{(2k)!},Q(X,Y)^k\Big\rangle
 d\nu(x,y)\\
 &&\textrm{(by the equality (\ref{pairing}) and $Q \in \mathcal{F}_{Np}$)}\\
 &=& \chi_{N}(Q)\chi_{p}(Q) \cdot \Big\langle \kappa(r)\cdot \int_{\bZ_p^{\times}\times\bZ_p}
 \chi_p^2(x)\cdot \frac{(xY-yX)^{2k}}{(2k)!}d\nu(x,y), Q^k \Big\rangle \\
 &=& \chi(Q) \cdot \langle \phi_{\kappa}(\mu), Q^k\rangle.
\end{eqnarray*}
This completes the proof.

\end{proof}

Now we are ready to define the overconvergent Shintani lifting
following the definition in \cite{St1}. 

\begin{defn}
For each $\Phi \in H^1_{c}(\Gamma_{0}(Np),
\cD_{N})=\mathrm{Hom}_{\Gamma_{0}(Np)}(\Delta_{0},\cD_{N})$ and
each $Q \in \cF_{Np}$ we define
\begin{eqnarray}
J(\Phi, Q):=\tilde{J}_Q(\Phi(D_Q)) \in \DN
\end{eqnarray}
where $\Phi(D_Q)$ is the value of the cocycle representing $\Phi$
on $D_Q$. \label{def16}
\end{defn}
Note that this definition of $J(\Phi,Q)$ does not depend on the
representative cocycle. A simple calculation confirms that
\begin{eqnarray*}
\tilde{J}_Q\big(\Phi(D_Q)\big)=\tilde{J}_{Q\mid
\gamma}\big(\Phi(D_{Q\mid \gamma})\big),
\end{eqnarray*}
for $\gamma\in \Gam$ and $D_Q=\partial C_Q$ as in Definition
\ref{cycle}. Thus $J(\Phi, Q)$ depends only on the
$\Gamma_{0}(Np)$-equivalence class of $Q$ and the following
definition makes sense.
\begin{defn}
We define the overconvergent Shintani lifting
$\Theta:H^1_{c}(\Gamma_{0}(Np), \cD_{N}) \to \DN[[q]]$ by

\begin{displaymath}
\Theta(\Phi):=\left\{ \begin{array}{ll} \sum_{Q \in
\mathcal{F}_{Np}/\Gamma_{0}(Np)}J(\Phi,Q)q^{\delta_{Q}/Np} & \textrm{if $Np$ is odd,}\\
\\
\sum_{Q \in
\mathcal{F}_{Np}/\Gamma_{0}(Np)}J(\Phi,Q)q^{\delta_{Q}/4Np} &
\textrm{if $Np$ is even.}
\end{array} \right.\label{maindef}
\end{displaymath}

\end{defn}

We have seen that $\Hc$ is an $\cH$-module but $\DN[[q]]$ doesn't
have an $\cH$-action a priori. We will put a $\cH$-module
structure on $\DN[[q]]$ so that $\Theta$ is $\cH$-equivariant. In
order to do so, we need an explicit description of
$\Theta(\Phi\big|T_l)$ for all primes $l$ and
$\Theta(\Phi\big|T_{l,l})$ for all primes $l\nmid Np$ in terms of
the coefficients of the $q$-expansion of $\Theta(\Phi)$ (see the
(\ref{Tpp}) for the definition of $T_{l,l}$) since $\cH$ is
generated by these Hecke operators over $\NS$. For $\Phi \in \Hcb$
and $m \in \Delta_N$, we define the notation $\Phi \otimes m \in
\Hc$ as follows:
\begin{eqnarray*}
(\Phi \otimes m)(D)=\Phi(D) \ten m \in \cD \ten \NS=:\cD_N
\end{eqnarray*}
for any $D \in \Delta_0$. Even though we state the theorem for
$\Phi \otimes 1$ for simplicity, the formula for general $\Phi \in
\Hc$ is also straightforward as $\Theta$ is $\NS$-linear map.

\begin{thm}
First, $\Theta$ is a $\NS$-linear map. Also, if $\tilde{\Phi}=\Phi
\otimes 1 \in \Hc$ is an overconvergent modular symbol and we
write
\begin{eqnarray*}
\Theta(\tilde{\Phi})=\sum_{n=1}^{\infty}\big(1\TNS a(n)\big)q^n
\in \DN[[q]],
\end{eqnarray*}
 then we have the following explicit formula for $T_l$ on the $q$-expansion of
 $\Theta(\tilde{\Phi})$:
 \begin{eqnarray}
 \Theta(\tilde{\Phi}\big|T_l)=
 \sum_{n=1}^{\infty}\Big(1\TNS\big(a(nl^2)+(\frac{Np\cdot
 n}{l})[l]_N\delta_l \ast a(n)+l[l^2]_N\delta_{l^2}\ast
 a(\frac{n}{l^2})\big)\Big)q^n
 \label{Hecke}
 \end{eqnarray}
 for any odd prime $l$. Here $\delta_s\in \cD(\bZ_p^{\times})$ is the Dirac distribution
 at $s \in \bN$ (we put $\delta_s=0$, if $p$ divides $s$), $*$ is the convolution product
 in $\cD(\bZ_p^{\times})$, and we put $a(\frac{n}{l^2})=0$ for $n$ not divisible by $l^2$.
 We also have the formula for $T_{l,l}$ for any prime number $l \nmid
 Np$
 \begin{eqnarray}
 \Theta(\tilde{\Phi}\big|T_{l,l})=
 \sum_{n=1}^{\infty}\big(1 \TNS [l^2]_N\cdot\delta_{l^2}\ast
 a(n)\big)q^n
 \label{Hecke1}
 \end{eqnarray}with the same notation as above.
\label{main1}
\end{thm}

\begin{proof}
The $\NS$-linearity of $\Theta$ follows from Proposition \ref{JQ}.
We first concentrate on the explicit description of
$\Theta(\tilde{\Phi}\big|T_l)$. The description of $T_{l,l}$ for
$l\nmid Np$ will be much easier. For $Np$ odd (resp. even) denote
by $d_n$ the discriminant of the number field $K=\bQ(\sqrt{Np\cdot
n})$ (resp. $K=\bQ(\sqrt{4Np\cdot n}))$ and put $Np\cdot
n=d_nc_n^2$ (resp. $4Np\cdot n=d_nc_n^2)$. The positive number
$c_n$ defined by the above equations is a positive integer or a
positive half integer but we will consider only the positive
integer $c_n$ (if $c_n$ is not an integer, the terms appearing
below are 0). We let

\begin{displaymath}
\cL_s(Np) = \left\{ \begin{array}{ll} \{Q=aX^2+bXY+cY^2\in
\cF_{Np} \ \big|\
(a,b,c)=1,\ \delta_Q=Np\cdot s\} & \textrm{if $Np$ is odd}\\
\\
\{Q=aX^2+bXY+cY^2\in \cF_{Np} \ \big| \ (a,b,c)=1,\
\delta_Q=4Np\cdot s\} & \textrm{if $Np$ is even}
\end{array} \right.
\end{displaymath}for $s \in \bN$. Then, from the construction of $\Theta$, we have

\begin{displaymath}
\Theta(\tilde{\Phi})=\sum_{n=1}^{\infty}\Bigg(\sum_{\substack{m|c_n,m>0\\
(m,Np)=1 }}\bigg(\sum_{\substack{Q\ \textrm{belongs to} \\
\cL_{\frac{n}{m^2}}(Np)/\Gam}}\Big(1\TNS
J_{mQ}(\Phi(D_{Q}))\cdot[mQ]_N\Big)\bigg)\Bigg)q^n.
\end{displaymath}
Notice that $D_{Q}=D_{mQ}$ for $m \in \bN$ satisfying $(m,Np)=1$
(see the Lemma 2.7. (iii) of \cite{Shn}). Put
\begin{displaymath}
\alpha_{_{Np}}(s,\tilde{\Phi},m):=\sum_{\substack{Q\ \textrm{belongs to} \\
\cL_s(Np)/\Gam}}\bigg(J_{mQ}(\Phi(D_{Q}))\cdot[mQ]_N\bigg).
\end{displaymath}
In order to get the desired description of $T_l$ on
$\Theta(\tilde{\Phi} )$, we will first derive a formula for
$\alpha_{_{Np}}(s,\tilde{\Phi}\big|T_l,m)$.
\begin{lem}
Notations being as above, for $s\in \bN$, we have
\begin{displaymath}
\alpha_{_{Np}}(s,\tilde{\Phi}\big|T_l,m)=\left\{ \begin{array}{ll}
\alpha_{_{Np}}(sl^2,\tilde{\Phi},m)+\big(1+(\frac{d_s}{l})\big)[l]_N \delta_l*\alpha_{_{Np}}(s,\tilde{\Phi},m) & \textrm{if}\ (l,c_s)=1\\
\\
\alpha_{_{Np}}(sl^2,\tilde{\Phi},m)+\big(l-(\frac{d_s}{l})\big)[l^2]_N
\delta_{l^2}*\alpha_{_{Np}}(\frac{s}{l^2},\tilde{\Phi},m) &
\textrm{if}\ (l^2,c_s)=l \\
\\
\alpha_{_{Np}}(sl^2,\tilde{\Phi},m)+l[l^2]_N
\delta_{l^2}*\alpha_{_{Np}}(\frac{s}{l^2},\tilde{\Phi},m) &
\textrm{if}\ (l^2,c_s)=l^2
\end{array} \right.
\end{displaymath}
where $\delta_l$ and $\delta_{l^2}$ are defined in the theorem
\ref{main1} and we put
$\alpha_{_{Np}}(\frac{s}{l^2},\tilde{\Phi},m)=0$ for $s$ not
divisible by $l^2$. \label{main lem}
\end{lem}
\begin{proof}
Let
\begin{eqnarray*}
\alpha_0=\left(
           \begin{array}{cc}
             l & 0 \\
             0 & 1 \\
           \end{array}
         \right)\
         \textrm{and}\ \alpha_j=\left(
                                  \begin{array}{cc}
                                    1 & j \\
                                    0 & l \\
                                  \end{array}
                                \right)
                                \textrm{for $1\leq j \leq l$}
\end{eqnarray*}

Then, for $T_l \in \cH$,
\begin{displaymath}
\tilde{\Phi}\big|T_l=\sum_{j=0}^{l}\tilde{\Phi}|\alpha_j.
\end{displaymath}
Having in mind the action of $\alpha_j$ on $\cD \ten \NS=\cD_N$,
we have
\begin{displaymath}
\alpha_{_{Np}}(s,\tilde{\Phi}\big|T_l,m)=\sum_{\substack{Q\
\textrm{belongs to} \\  \cL_s(Np)/\Gam}}\sum_{j=0}^{l} J_{mQ|
\alpha_j^{-1}}(\Phi( \alpha_j\cdot
D_Q))\cdot[mQ]_N[a(\alpha_j)^2]_N,
\end{displaymath}
where $a(\alpha_j)$ is the upper left entry of $\alpha_j$.


Assume $d_s>1$ and we have $K=\bQ(\sqrt{d_s})=\bQ(\sqrt{d_s
c_s^2})$ (remember that $d_s$ is the discriminant of $K$). In this
case $Q\in \cL_s(Np)$ has a non-perfect square discriminant and so
$D_Q=\{\gamma_Q \omega\}-\{\omega\}$ for any point $\omega \in
\mathfrak{h}$ by the Lemma 2.7 (i) of \cite{Shn}, where $\gamma_Q$
was given in the Definition \ref{cycle}. Let $Q_1,Q_2,...,Q_h$ be
a complete set of representatives of $\Gam$-equivalence classes of
quadratic forms in $\cL_s(Np)$. Let
\begin{displaymath}
Q_i(X,Y)=a^iX^2+b^iXY+c^iY^2 \ (1\leq i \leq h).
\end{displaymath}
Also put $\omega_1^i=\frac{b^i+c_s \sqrt{d_s}}{2}$ and
$\omega_2^i=c^i$ and denote by $L^i$ the lattice in $K$ generated
by $\omega_1^i$ and $\omega_2^i$. If we denote by
$\cL^i=\cL^i(X,Y)$ the primitive integral binary quadratic form
with discriminant $d_sc_s^2$ associated to $L^i$ (look at the
Lemma 2.2 of \cite{Shn} for the precise dictionary between the
lattices in $K$ and quadratic forms.). In fact it's easy to see
\begin{eqnarray*}
\cL^i(X,Y)=Q_i(X,Y) \ (1 \leq i \leq h).
\end{eqnarray*}
It follows from the Lemma 2.2 and Lemma 2.6 of \cite{Shn} that
$\{L^1,L^2,...,L^h\}$ forms a complete set of representatives of
equivalence classes of lattices in $K$ with conductor $c_s$ (see
the page 110 of \cite{Shn} for the defintion of the equivalence
relation and the conductor of $L^i$).\\
Put
\begin{eqnarray*}
Q_{i,j}:=Q_i|\alpha_j^{-1} \ (1\leq i \leq h \ \textrm{and} \
0\leq
j \leq l) \\
 e_1=[\cO_{c_s}^{\times}:\cO_{c_sl}^{\times}] \
\textrm{and} \ e_2=[\cO_{c_s}^{\times}/(c_s,l):\cO_{c_s}^{\times}]
\end{eqnarray*}
where the definition of $\cO_{c_s}^{\times}$ is given on the page
111 of \cite{Shn}. Having in mind that
\begin{displaymath}
[Q_{i,j}]_{N}=[Q_i]_{N} \ \ (1 \leq j \leq l) \ \textrm{and} \
[Q_{i,j}]_{N}=[l^2]_{N}\cdot [Q_i]_{N} \  (j=0),
\end{displaymath}
we get
\begin{eqnarray*}
\alpha_{_{Np}}(s,\tilde{\Phi}\big|T_l,m)&=&\sum_{i=0}^{h}\sum_{j=0}^{l}J_{mQ_{i,j}}\bigg(\Phi
\big(\{\alpha_j \gamma_{Q_i} \alpha_j^{-1}\cdot
\omega\}-\{\omega\}\big)\bigg)\cdot[mQ_{i,j}]_N \\
&=&\sum_{i=0}^{h}e_1^{-1}\sum_{j=0}^{l}J_{mQ_{i,j}}\bigg(\Phi
\big(\{\alpha_j \gamma_{Q_i}^{e_1} \alpha_j^{-1}\cdot
\omega\}-\{\omega\}\big)\bigg)\cdot[mQ_{i,j}]_N.
\end{eqnarray*}
Define
\begin{displaymath}
L^{i,j}:=\left\{ \begin{array}{ll}
l \omega^i_1\bZ +\omega_2^i\bZ & \textrm{if $j=0$}\\
\\
(\omega_1^i + j \omega_2^i)\bZ+l \omega_2^i\bZ & \textrm{if $
1\leq j \leq l$}.
\end{array} \right.
\end{displaymath}
Then $L^{i,0},L^{i,1},...,L^{i,l}$ are mutually distinct
sublattices of $L^i$ with index $l$. Let $\cL^{i,j}$ be the
quadratic form corresponding to $L^{i,j}$. Then we have
\begin{displaymath}
Q_{i,j}(X,Y)=\left\{ \begin{array}{ll}
\cL^{i,j}(X,Y) & \textrm{if the conductor of $L^{i,j}=c_sl$}\\
\\
l\cL^{i,j}(X,Y) & \textrm{if the conductor of $L^{i,j}=c_s$}\\
\\
l^2\cL^{i,j}(X,Y) & \textrm{if the conductor of
$L^{i,j}=\frac{c_s}{l^2}$}
\end{array} \right.
\end{displaymath}
and
\begin{displaymath}
\alpha_j \gamma_{Q_i}^{e_1}\alpha_j^{-1}=\left\{ \begin{array}{ll}
\gamma_{\cL^{i,j}} & \textrm{if the conductor of $L^{i,j}=c_s l$}\\
\\
\gamma_{\cL^{i,j}}^{e_1} & \textrm{if the conductor of $L^{i,j}=c_s$}\\
\\
\gamma_{\cL^{i,j}}^{e_1e_2} & \textrm{if the conductor of
$L^{i,j}=\frac{c_s}{l^2}$}
\end{array} \right.
\end{displaymath}

Now our Lemma \ref{main lem} can be derived from Lemma 2.5, Lemma
2.7 and Lemma 2.2 of \cite{Shn}. For example, if $(l,c_s)=1$, then
Lemma 2.5, Lemma 2.7 and Lemma 2.2 of \cite{Shn} implies the
following:

\begin{eqnarray*}
&{}&\alpha_{_{Np}}(s,\tilde{\Phi}\big|T_l,m)\\
&=&e_1e_1^{-1}\sum_{\substack{\cL^{i,j}\
\textrm{belongs to}\\
\cL_{sl^2}(Np)/\Gam}}J_{m\cL^{i,j}}(\Phi(D_{\cL^{i,j}}))[m\cL^{i,j}]_N
+\big(1+(\frac{d_s}{l})\big)\sum_{\substack{\cL^{i,j}\
\textrm{belongs to} \\
\cL_{s}(Np)/\Gam}}J_{lm\cL^{i,j}}(\Phi(D_{i,j}))[lm\cL^{i,j}]_N \\
&=&\alpha_{_{Np}}(sl^2,\tilde{\Phi},m)+\big(1+(\frac{d_s}{l})\big)[l]_N
\delta_l*\alpha_{_{Np}}(s,\tilde{\Phi},m).
\end{eqnarray*}

The cases $(l^2,c_s)=l$ and $(l^2,c_s)=l^2$ also can be checked in
the similar way. If $d_s=1$, the proof is similar and much
simpler. This completes the proof of Lemma \ref{main lem}.
\end{proof}

Now we will use the Lemma \ref{main lem} to prove the theorem. We
can express $\Theta(\tilde{\Phi})$ as follows:

\begin{displaymath}
\Theta(\tilde{\Phi})=\sum_{n=1}^{\infty}\big(1\TNS a(n)\big)q^n=\sum_{n=1}^{\infty}\Bigg(1\TNS\Big(\sum_{\substack{m|c_n,m>0\\
(m,Np)=1 }}\alpha_{_{Np}}(\frac{n}{m^2},\tilde{\Phi},m)\Big)
\Bigg)q^n.
\end{displaymath}
If we apply $T_l$, then we have

\begin{displaymath}
\Theta(\tilde{\Phi}\big|T_l)=\sum_{n=1}^{\infty}\Bigg(1\TNS\Big(\sum_{\substack{m|c_n,m>0\\
(m,Np)=1
}}\alpha_{_{Np}}(\frac{n}{m^2},\tilde{\Phi}\big|T_l,m)\Big)
\Bigg)q^n.
\end{displaymath}
Let
\begin{displaymath}
b(n):=\sum_{\substack{m|c_n,m>0\\
(m,Np)=1 }}\alpha_{_{Np}}(\frac{n}{m^2},\tilde{\Phi}\big|T_l,m).
\end{displaymath}
Note that we can assume $c_n$ is an integer. Let's analyze the
$q^n$-coefficient case by case. If $l$ is a factor of $Np$, we
have, by Lemma \ref{main lem},

\begin{eqnarray*}
b(n)&=&\sum_{\substack{m|c_n,m>0\\
(m,Np)=1 }}\alpha_{_{Np}}(\frac{nl^2}{m^2},\tilde{\Phi},m) \\
&=&\sum_{\substack{m|c_{nl^2},m>0\\
(m,Np)=1 }}\alpha_{_{Np}}(\frac{nl^2}{m^2},\tilde{\Phi},m)\\
&=&a(nl^2) ,
\end{eqnarray*}
since $l|p$ implies $\delta_l=0$ $(\delta_{l^2}=0$) and $l|N$
implies $[l]_N=0$ $([l^2]_N=0)$. Next we examine the case that $l$
is prime to both $Np$ and $c_n$. Then Lemma \ref{main lem} implies
that
\begin{eqnarray*}
b(n)&=&\sum_{\substack{m|c_{n},m>0\\
(m,Np)=1}}\big(\alpha_{_{Np}}(\frac{nl^2}{m^2},\tilde{\Phi},m)+\big(1+(\frac{d_{\frac{n}{m^2}}}{l})\big)[l]_N
\delta_l*\alpha_{_{Np}}(\frac{n}{m^2},\tilde{\Phi},m)\big)
\\
&=&\sum_{\substack{m|c_{n},m>0\\
(m,Np)=1 }}\alpha_{_{Np}}(\frac{nl^2}{m^2},\tilde{\Phi},m)+\sum_{\substack{m|c_{n},m>0\\
(m,Np)=1 }}\big(1+(\frac{d_{n}}{l})\big)[l]_N
\delta_l*\alpha_{_{Np}}(\frac{n}{m^2},\tilde{\Phi},m) \\
&=&\sum_{\substack{m|c_{nl^2}=lc_n,m>0\\
(m,Np)=1 }}\alpha_{_{Np}}(\frac{nl^2}{m^2},\tilde{\Phi},m)+\sum_{\substack{m|c_{n},m>0\\
(m,Np)=1 }}\big(\frac{Np\cdot n}{l}\big)[l]_N
\delta_l*\alpha_{_{Np}}(\frac{n}{m^2},\tilde{\Phi},m)\\
&=&a(nl^2)+\big(\frac{Np\cdot n}{l}\big)[l]_N \delta_l*a(n).
\end{eqnarray*}
Finally consider the case that $l$ is prime to $Np$ but is a
factor of $c_n$. Set $c_n=c_n'l^k$, where $c_n'$ is prime to $l$.
Then we get
\begin{eqnarray*}
b(n)&=&\sum_{\substack{m'|c_{n}',m'>0\\
(m',Np)=1 \\0\leq i \leq k \\ \textrm{or}\
m'=2,i=0}}\alpha_{_{Np}}(\frac{n}{m'^{2}l^{2i}},\tilde{\Phi}|T_l,m'l^i).
\end{eqnarray*}
Since Lemma \ref{main lem} shows that
\begin{displaymath}
\alpha_{_{Np}}(\frac{n}{m'^{2}l^{2i}},\tilde{\Phi}|T_l,t) =
\left\{
\begin{array}{ll}
\alpha_{_{Np}}(\frac{nl^2}{m'^{2}l^{2i}},\tilde{\Phi},t)+l[l^2]_N
\delta_{l^2}*\alpha_{_{Np}}(\frac{n}{m'^2l^{2i+2}},\tilde{\Phi},t)
& (0\leq i \leq k-2) \\
\\
\alpha_{_{Np}}(\frac{nl^2}{m'^{2}l^{2i}},\tilde{\Phi},t)+\big(l-(\frac{d_n}{l})\big)[l^2]_N
\delta_{l^2}*\alpha_{_{Np}}(\frac{n}{m'^2l^{2i+2}},\tilde{\Phi},t) & (i=k-1)\\
\\
\alpha_{_{Np}}(\frac{nl^2}{m'^{2}l^{2i}},\tilde{\Phi},t)+\big(1+(\frac{d_n}{l})\big)[l]_N
\delta_{l}*\alpha_{_{Np}}(\frac{n}{m'^2l^{2i}},\tilde{\Phi},t) &
(i=k)
\end{array} \right.
\end{displaymath}
for any $t\in \bN$, we can conclude that
\begin{eqnarray*}
b(n)&=&\sum_{\substack{m|c_{nl^2}=lc_n,m>0\\
(m,Np)=1 }}\alpha_{_{Np}}(\frac{nl^2}{m^2},\tilde{\Phi},m)+\sum_{\substack{m|c_{\frac{n}{l^2}}=\frac{c_n}{l},m>0\\
(m,Np)=1 }}l[l^2]_N
\delta_{l^2}*\alpha_{_{Np}}(\frac{n}{m^2l^2},\tilde{\Phi},m)\\
&=&a(nl^2)+l[l^2]_N\delta_{l^2}\ast
 a(\frac{n}{l^2}).
\end{eqnarray*}
This completes the proof of the description of $T_l$-action in the
theorem \ref{main1}. The formula $T_{l,l}$-action for $l \nmid Np$
is easily derived, because $T_{l,l}$-action on $\tilde{\Phi}$ is
the same as $\left(
               \begin{array}{cc}
                 l & 0 \\
                 0 & l \\
               \end{array}
             \right)$-action on it.

\end{proof}

Now we define the action of $T_{l}$ for any rational prime $l$ and
$T_{l,l}$ for $l\nmid Np$ on formal $q$-expansions
\begin{eqnarray*}
\Theta=\sum_{n=1}^{\infty}\big(a(n)\TNS b(n)\big)q^n \in \DN[[q]]
\end{eqnarray*}

as follows:
\begin{eqnarray}
\Theta\big|T_l=
 \sum_{n=1}^{\infty}\Big(a(n) \TNS\big(b(nl^2)+(\frac{Np\cdot
 n}{l})[l]_N\delta_l \ast b(n)+l[l^2]_N\delta_{l^2}\ast
 b(\frac{n}{l^2})\big)\Big)q^n
 \label{Hecke action1}
\end{eqnarray}
\begin{eqnarray}
\Theta\big|T_{l,l}=\sum_{n=1}^{\infty}\Big(a(n) \TNS [l^2]_N \cdot
\delta_{l^2}\ast
 b(n)\Big)q^n.
 \label{Hecke action2}
\end{eqnarray}
Then we can conclude that the overconvergent Shintani lifting
$\Theta$ is $\cH$-linear map by Theorem \ref{main1} except for
operators containing $T_2$. But if $N$ is even, then even
$\Theta(\Phi\big|T_2)=\Theta(\Phi)\big|T_2$ is true so that
$\Theta$ is actually $\cH$-linear.

Because we will show that our overconvergent Shintani lifting
$\Theta$ interpolates the classical Shintani liftings, it is
reasonable to define the image of $\Theta$ to be a universal
overconvergent half-integral weight modular form.
\begin{defn}
We say that $\Theta \in \DN[[q]]$ is a universal overconvergent
half-integral weight modular form if
\begin{eqnarray*}
\Theta=\Theta(\Phi)
\end{eqnarray*}
for some $\Phi \in \Hc$.\label{univmod}
\end{defn}
A universal overconvergent half-integral weight modular form will
be viewed as a $p$-adic family of overconvergent half-integral
weight modular forms. Finally we give a definition of an
overconvergent half-integral weight modular form of weight
$\kappa\in \textrm{Hom}_{cts}(\bZ_{N,p}^{\times},K^{\times})$.
\begin{defn}
We say that $\theta \in \bC_p[[q]]$ is an overconvergent
half-integral weight modular form of weight $\kappa\in
\textrm{Hom}_{cts}(\bZ_{N,p}^{\times},K^{\times})$ if it can be
written as
\begin{eqnarray*}
\theta=\Theta(\kappa)=\sum_{n=1}^{\infty}\big(1\TNS\alpha_n(\kappa)\big)q^n
\end{eqnarray*}
where $\Theta=\sum_{n=1}^{\infty}\big(1\TNS\alpha_n\big) q^n \in
\DN[[q]]$ is a universal overconvergent half-integral weight
modular form.\label{overmod}
\end{defn}

\section{Connection to the eigencurve and $p$-adic family of half-integral weight modular
forms}\

Since points on the (half-integral) eigencurve correspond to
(half-integral) overconvergent Hecke eigenforms, the existence of
a Hecke eigensymbol in $\Hc$ which interpolates overconvergent
Hecke eigenforms and a local version of overconvergent Shintani
lifting will give us a local piece of the half-integral eigencurve
(see \cite{Ram2}). So we will construct a Hecke eigensymbol and a
local version of lifting whose meaning will be precise later.

\subsection{Slope $\leq h$ decompositions of overconvergent modular symbols}
\

The goal of this section is to review the result of G. Stevens in
\cite{St3} which guarantees the existence of slope $\leq h$
decomposition of overconvergent modular symbols over
some $K$-affinoid space.\\

In order to make a connection to the (half-integral) eigencurve,
we need an appropriate $p$-adic Hecke algebra over
$\bZ_p[[\bZ_{p,N}^{\times}]]$, but unlike the ordinary case (slope
0 case), it is impossible to define a global $p$-adic Hecke
algebra finite over $\bZ_p[[\bZ_{p,N}^{\times}]]$ which
parametrizes all overconvergent Hecke eigenforms of arbitrary
finite slope. Instead we can define a local $p$-adic Hecke algebra
finite over $A(B)$, where $B$ is some affioid subdomain of the
weight space $\cX_N$, using slope $\leq h$ decomposition. We first
review some background knowledge on rigid spaces, especially
weight
space, and state the theorem of G. Stevens in \cite{St3}. \\

Let $\Lambda$ be the completed group algebra
$\bZ_p[[\bZ_p^{\times}]]$ and let $\Lambda_N$ be the completed
group algebra $\bZ_p[[\bZ_{p,N}^{\times}]]$ where
$\bZ_{p,N}^{\times}:=\bZ_p^{\times}\times(\bZ/N\bZ)^{\times}=\bZ_p^{\times}\times
\Delta_N$. Then $\bZ_p[[\bZ_{p,N}^{\times}]] \cong
\bZ_p[[\bZ_p^{\times}]][(\bZ/N\bZ)^{\times}]\cong
\bigoplus_{\chi}\bZ_p[[\bZ_p^{\times}]]$ where the direct product
is taken over the Dirichlet characters of $(\bZ/N\bZ)^{\times}$.

Now consider the category $AF_{\bZ_p}$ whose objects are $Spf(A)$
where $A$ is an adic and noetherian ring  which is $\bZ_p$-algebra
and $A/J$ ($J$ is the biggest ideal of definition) is a finitely
generated $\mathbb{F}_p$-algebra and morphisms are locally ringed
space morphisms. Let $Rig_{K}$ be the category of $K$-rigid
analytic varieties.($K$-rigid analytic variety is a locally
$G$-ringed space and a morphism is a locally $G$-ringed morphism.-
A locally $G$-ringed space is a pair $(X,\mathcal{O}_X)$, where
$X$ is a set equipped with a saturated Grothendieck topology and
$\mathcal{O}_X$ is a sheaf of rings such that all stalks are local
rings.) P.Berthelot had constructed a functor from $AF_{\bZ_p}$ to
$Rig_K$. We refer to [Col-Mz] and [De] for this construction. Note
that $Spf(\bZ_p[[1+p\bZ_p]])$ belongs to $AF_{\bZ_p}$. So we can
construct a $K$-rigid analytic variety associated to
$\bZ_p[[1+p\bZ_p]]$. This turns out to be the open unit disk
$B(0,1)_K$ defined over $K$. Since $\Lambda$ and $\Lambda_N$ are
finite direct sums of $\bZ_p[[1+p\bZ_p]]$, we can apply P.
Bertholet's contruction to $\Lambda$ and $\Lambda_N$. We denote
the resulting $K$-rigid analytic varieties by $\cX$ and $\cX_N$
respectively. In this case $\cX$ (respectively $\cX_N$) is the
finite union of $\varphi(p)$ (respectively $\varphi(pN)$) open
unit disks and each component corresponds to the Dirichlet
character of $(\bZ/p\bZ)^{\times}$ (respectively
$(\bZ/pN\bZ)^{\times}$). We call $\cX$ (respectively $\cX_N$)
\textbf{weight space} (respectively \textbf{weight space}
of tame conductor $N$). So $\cX=\cX_1$. \\

For any $K$-rigid analytic variety $X$ (respectively any
$K$-affinoid variety $X$) we define $A(X)$ to be the ring of
$K$-rigid analytic functions on $X$ (respectively $K$-affinoid
functions on $X$). Then there is a natural ring homomorphism
\begin{eqnarray}
 \Lambda_N \to A(\cX_N)
\label{natural}
\end{eqnarray}
for $N \geq 1$. In the case that $X$ is a $K$-affinoid variety
over $\cX_N$, the $R$-valued points of $X$, for any commutative
$\bQ_p$-algebra $R$, are given by continuous homomorphisms from
$A(X)$ to $R$, i.e.
\begin{eqnarray*}
X(R):=\textrm{Hom}_{cts}(A(X),R).
\end{eqnarray*}
The $R$-valued points of $\cX_N$ are given by
\begin{eqnarray*}
\cX_N(R)=\textrm{Hom}_{cts}(\Lambda_N,R).
\end{eqnarray*}
for $N \geq 1$. Recall that there is a natural bijection between
\begin{eqnarray}
\textrm{Hom}_{cts}(\Lambda_N,\bC_p) &\textrm{and}&
\textrm{Hom}_{cts}(\bZ_{p,N}^{\times},\bC_p^{\times})
\label{points}.
\end{eqnarray}
So in particular, the $\bC_p$-points of the weight space $\cX_N$
is $\textrm{Hom}_{cts}(\bZ_{p,N}^{\times},\bC_p^{\times})$.

\begin{defn}
A character $\kappa:1+p\bZ_p \to \bC_p^{\times}$ is called
arithmetic if
\begin{eqnarray*}
\kappa(t)=t^k
\end{eqnarray*}
for some $k\in \bZ^{\geq0}$ and for all $t$ sufficiently close to
1 in $1+p\bZ_p$. \label{def2}
\end{defn}
\begin{defn}
Let $\kappa \in X(\bC_p)$ be a $\bC_p$-point of a $K$-affinoid
variety $X$ over $\cX_N$ or $\cX_N$ itself. We call $\kappa$ an
arithmetic point of signature of $(k,\chi)$, if the associated
group character $\kappa:\bZ_{p,N}^{\times} \to \bC_p^{\times}$,
via the natural homomorphism (\ref{natural}) and the
correspondence (\ref{points}), satisfies
\begin{eqnarray*}
\kappa(t)=\chi(t)\cdot t_p^k
\end{eqnarray*}
where $t=(t_p,t_N)\in \bZ_{p}^{\times}\times(\bZ/N\bZ)^{\times}$
and for some $\chi$, a finite character of $\bZ_{p,N}^{\times}$,
and $k\in \bZ^{\geq0}$. \label{def3}
\end{defn}
We use the following notation: for any $K$-affinoid variety $X$
(or $\cX_N$) and any commutative $\bQ_p$-algebra $R$,
\begin{eqnarray*}
X^{arith}(R):=\{\kappa\in X(R) \ | \ \kappa \ \textrm{is
arithmetic.} \}.\label{arith notation}
\end{eqnarray*}
Now we state the connection between locally analytic distributions
$\NS$ and $A(\cX_N)$, the $K$-rigid analytic functions on $\cX_N$.
Recall that there is a topological $K$-algebra isomorphism (as
$K$-Fr\'{e}chet spaces)
\begin{eqnarray}
F:\cD(\bZ_{p,N}^{\times}) &\simeq& A(\cX_N) \label{Fourier} \\
\mu \ \ & \mapsto & \ f_{\mu}
\end{eqnarray}where $f_{\mu}(\kappa):=\int_{\bZ_{p,N}^{\times}}\kappa(x)
d\mu(x)$ which is the $p$-adic Fourier transform of Amice (see the
\cite{Ami1} and \cite{Ami2}).

\begin{defn}
We define
\begin{equation}
\cD_{B}:=\cD \otimes_{A(\cX)} A(B)
\end{equation}
where $B$ is any $K$-affinoid subdomain of $\cX_N$ and $A(B)$ is
the associated $K$-affinoid algebra, endowed with the spectral
norm $\mid\cdot\mid_{A(B)}$. \label{def4}
\end{defn}

Let $S_0(Np)$ act on $\cD_{B}$ by the formula $(\mu
\otimes_{A(\cX)}\lambda)|\gamma:=\mu|\gamma\otimes_{A(\cX)}[a]_{N}\lambda$
for
$B\subseteq \cX_N$ where $\gamma=\left(%
\begin{array}{cc}
  a & b \\
  c & d \\
\end{array}%
\right)$ and $[a]_{N}$ is the image of $a$ in $
(\bZ/N\bZ)^{\times}\subseteq \Lambda_N^{\times}$. Since the action
of $A(B)$ on $\cD_B$ commutes with the action of $\Gam$ on
$\cD_B$, $\HB$ inherits an action of $A(B)$. Indeed $\HB$ is
equipped with an $A(B)[T_n: n\in\bN]$-module structure. In
particular, $U_p$ acts on $H^1_c(\Gamma_0(N), \cD_B)$. In fact we
have an $A(B)$-linear endomorphism $U_p:\HB \to \HB$. G. Stevens
proved the following theorem in [St3] which is essential to prove
the existence of slope $\leq h$ decompostion.

\begin{thm}
\
\begin{quote}
(1) $H^1_c(\Gamma_0(Np), \cD_B)$ is orthonormalizable as an
$A(B)$-module.
\end{quote}
\begin{quote}
(2) The $U_p$-action on $H^1_c(\Gamma_0(Np), \cD_B)$ is completely
continuous.
\end{quote}\label{thm1}
\end{thm}

The above theorem enables us to apply the Fredholm-Riesz-Serre
theory, and in particular, we have the Fredholm series
$\det(1-U_p\cdot T|\HB) \in A(B)[[T]]$ which is entire in $T$. For
details
see \cite{St3} and \cite{Ser}.\\

Now we recall the definition of slope $\leq h$ decomposition over
$B$ from [A-S2].
\begin{defn}
$x \in \HB$ is said to be have slope $\leq h$ with respect to
$U_p$ for some $h \in \bR$, if there is a polynomial $Q \in
A(B)[T]$ with the following properties:\

$(1) Q^{\ast}(U_p)\cdot x=0$ where $Q^{\ast}(T):=T^d\cdot Q(1/T)$
with $d=deg(Q)$;\

$(2)$ the leading coefficient of $Q$ is a multiplicative unit with
respect to $|\cdot|_{A(B)}$; and \

$(3)$ every slope of $Q$ is $\leq h$, where the slopes of $Q$ are
the slopes of the Newton polygon of $Q$ (see [A-S2] for the
definition of the Newton polygon and its slope). $\HB^{(\leq h)}$
is defined as the set of all elements of $\HB$ having slope $\leq
h$. \label{def5}
\end{defn}

Then $\HB^{(\leq h)}$ is an $A(B)$-submodule of $\HB$.
\begin{defn}
A slope $\leq h$ decomposition of $\HB$ over $B$ is an
$A(B)[U_p]$-module decomposition
\begin{eqnarray*}
\HB=\HB^{(\leq h)} \oplus \HB^{(\leq h)}_{\ast}
\end{eqnarray*}
such that\

$(1)$ $\HB^{(\leq h)}$ is a finitely generated $A(B)$-module; and\

$(2)$ for every polynomial $Q \in A(B)[T]$ of slope $\leq h$, the
map
\begin{eqnarray*}
Q^{\ast}(U_p):\HB^{(\leq h)}_{\ast} \to \HB^{(\leq h)}_{\ast}
\end{eqnarray*}
is an isomorphism of $A(B)$-modules.\label{def6}
\end{defn}

Then the following theorem is due to G. Stevens(\cite{St3}).
\begin{thm}
Let $x_0 \in \cX_N(K)$ be any $K$-point and let $h$ be a fixed
nonnegative rational number. Then there exists a $K$-affinoid
subdomain $\BN \subseteq \cX_N$ containing $x_0$ such that
$H^1_c(\Gamma_{0}(Np), \cD_{\BN})$ admits a slope $\leq h$
decomposition over $\BN$.\label{thmlocal}
\end{thm}

$\HBhN^{\leq h}$ will play a central role to define the $p$-adic
overconvergent Hecke algebra with slope $\leq h$ whose maximum
spectrum (associated $K$-affinoid variety) can be thought of as a
local piece of the eigencurve.

\subsection{Slope $\leq h$ overconvergent $p$-adic Hecke algebra and $p$-adic families of overconvergent
Hecke eigenforms}\

Note that the Hecke operators $T_n$ for all $n\in \bN$ preserves
$\HBhN^{(\leq h)}$.

\begin{defn}
We define the universal overconvergent $p$-adic Hecke algebra of
tame conductor $N$ with slope $\leq h$
\begin{eqnarray}
\cR_{h,N}:=Im(\cH_B {}^{H \atop \To}
 End_{A(B_{h,N})}(H^1_c(\Gamma_0(Np),\cD_{B_{h,N}})^{(\leq h)}))
\end{eqnarray}
where $\cH_B:=A(B_{h,N})[T_n: n\in \bN]$ is the abstract
overconvergent Hecke algebra over $\BN$ of tame conductor
$N$.\label{def8}
\end{defn}

It is known that $A(B_{h,N})$ is a principal ideal domain and we
have the following proposition.

\begin{prop} We have that
\begin{quote}
(1) $\cR_{h,N}$ is a finite, flat and torsion-free
$A(B_{h,N})$-algebra.\\
(2) $\cR_{h,N}$ is a $K$-affinoid algebra.
\end{quote}\label{prop3}
\end{prop}
\begin{proof}
\

(1) Since $H^1_c(\Gamma_0(Np),\cD_{B_h})^{(\leq h)}$ is a finitely
generated $A(B_{h,N})$-module, $\cR_{h,N}$ should be a finitely
generated $A(B_{h,N})$-module by the definition. The
$A(B_{h,N})$-torsion-freeness of $\cR_{h,N}$ is clear from the
definition. Since $A(B_{h,N})$ is a PID, torsion-freeness
implies the freeness. So in particular $\cR_{h,N}$ is flat over $A(B_{h,N})$. \\

(2) This part follows from [BGR](Proposition 5 on page 223) since,
by (1),  $\cR_{h,N}$ is finitely generated as an
$A(B_{h,N})$-module.
\end{proof}

Let $\Omega_h:=Sp(\cR_h)$ (respectively,
$\Omega_{h,N}:=Sp(\cR_{h,N})$) be the $K$-affinoid space
associated to $\cR_h$ (respectively, $\cR_{h,N}$) - See 7.1 of
\cite{BGR}. Then we have $A(\Omega_h)=\cR_h$ and
$A(\Omega_{h,N})=\cR_{h,N}$. So we get the following commutative
diagrams:
\begin{displaymath}
\xymatrix{ \cR_{h,N} & \cR_{h} \ar[l]   \\
A(B_{h,N}) \ar[u]  & A(B_{h}) \ar[l] \ar[u]    \\
A(\cX_{N}) \ar[u] & A(\cX) \ar[u] \ar[l]       \\
 \Lambda_N \ar[u] & \Lambda \ar[u] \ar[l] } \ \ \ \ \ \ \ \ \ \
\xymatrix{ \Omega_{h,N} \ar[d] \ar[r] & \Omega_{h} \ar[d]  \\
B_{h,N} \ar[d] \ar[r]  &  B_{h} \ar[d]\\
\cX_{N}  \ar[r]  & \cX    }
\end{displaymath}

\begin{defn}
We define the universal overconvergent Hecke eigenform of tame
conductor $N$ with slope $\leq h$ to be
\begin{eqnarray}
\mathbf{f}_{h,N}:=\sum_{n=1}^{\infty}\alpha_n q^n \in
\cR_{h,N}[[q]]
\end{eqnarray}
where $\alpha_n$ is the image of $T_n$ in
$End_{A(B_{h,N})}(H^1_c(\Gamma_0(Np),\cD_{B_{h,N}})^{(\leq h)})$
under the map $H$.\label{def9}
\end{defn}

We recall that $S_{k+2}(\Gamma_0(Np^m),\chi;R)$ is the subspace of
$S_{k+2}(\Gamma_0(Np^m),\chi)$ consisting of cusp forms whose
Fourier coefficients belong to $R$ where $R$ is any subfield of
$\bC_p$. The following theorem was proven by R. Coleman [Col2,
0.1].

\begin{thm}
Let $\kappa \in \Omega^{arith}_{h,N}(K)$ be an arithmetic
$K$-point of signature $(k,\chi)$ with $k \geq 0$. If $h<k+1$ or
if $h=k+1$, and
$\mathbf{f}_{h,N}(\kappa):=\sum_{n=1}^{\infty}\alpha_n(\kappa)
q^n$ is not in the image of the map $\mathbf{\Theta}^{k+1}$ where
$\mathbf{\Theta}$ is the operator which acts as $q\frac{d}{dq}$ on
formal $q$-expansions , then
\begin{eqnarray*}
\mathbf{f}_{h,N}(\kappa):=\sum_{n=1}^{\infty}\alpha_n(\kappa) q^n
\in M_{k+2}(\Gamma_0(Np^m),\chi;K)
\end{eqnarray*}
where $\alpha_n(\kappa):=\kappa(\alpha_n)$ and $\chi$ is defined
modulo $Np^m$ and moreover, $\mathbf{f}_{h,N}(\kappa)$ is an Hecke
eigenform. Furthermore, the $U_p$-eigen value of
$\mathbf{f}_{h,N}(\kappa)$ has $p$-adic valuation (called the
slope of $\mathbf{f}_{h,N}(\kappa)$) $\leq h$.\label{thm3}
\end{thm}

We will use the following notations:
\begin{eqnarray*}
\ONc(K):=\{\kappa:\cR_{h,N} \to K \ | \ \mathbf{f}_{h,N}(\kappa)\
\textrm{is a classical Hecke eigenform.}\}.
\end{eqnarray*}
\begin{eqnarray*}
B_{h,N}^{classical}(K):=\{\hat{\kappa}:A(B_{h,N}) \to K \ | \
\exists \kappa \in \ONc(K)\ \textrm{such that} \ \kappa \
\textrm{is lying above}\ \hat{\kappa} \}.
\end{eqnarray*}

\begin{displaymath}
\xymatrix{\cR_{h,N} \ar[r]^{\kappa} & K \\
A(B_{h,N}) \ar[u] \ar[ur]_{\hat{\kappa}}}
\end{displaymath}

The elements of $\ONc(K)$ will be referred to as classical
$K$-points of $\ONc$. We can view $\mathbf{f}_{h,N}$ as an
analytic function on $\Omega_{h,N}$ which interpolates the
$q$-expansions of Hecke eigenforms at classical points.

\subsection{Overconvergent Hecke eigensymbol}
\

In this section we will prove the existence of an overconvergent
Hecke eigensymbol.

\begin{defn}
For $\kappa \in \ONc(K)$ and $\Phi\in \HBhN$, we define
\begin{eqnarray}
\Phi_\kappa:=\phi_{\kappa,\ast}(\Phi)
\end{eqnarray}
where the map $\phi_{\kappa,\ast}$ is given in (\ref{speci}). For
$\kappa \in \ONc(K)$ of signature $(k,\chi)$, we also define
\begin{eqnarray}
\varphi_\kappa^{\pm}:=\varphi_{\mathbf{f}_{h,N}(\kappa)}^{\pm}:=
\frac{1}{\Omega^{\pm}_{\mathbf{f}_{h,N}(\kappa)}}\psi^{\pm}_{\mathbf{f}_{h,N}(\kappa)}
\end{eqnarray}
where we fix, once and for all, complex periods $\Omega_{f}^{\pm}$
for each Hecke eigenform $f$ so that
\begin{eqnarray}
\frac{1}{\Omega^{\pm}_{\mathbf{f}_{h,N}(\kappa)}}\cdot
\psi^{\pm}_{\mathbf{f}_{h,N}(\kappa)} \in
H^1_{par}(\Gamma_0(Np^m),L_{k,\chi}(R_\kappa))\label{algbraic}
\end{eqnarray}
For details see (4.3.4) of \cite{St1} about the existence of these
complex periods. \label{def10}
\end{defn}

Once we fix complex periods for each $f$, we can give an algebraic
version of the classical Shintani $\theta$-lifting. Define
\begin{eqnarray}
\theta_{k,\chi}^{\ast}(f):=\frac{1}{\Omega_f^-}\cdot\theta_{k,\chi}(f)\label{algShin}
\end{eqnarray}
for $f \in S_{2k+2}(\Gamma_0(M),\chi^2)$, a Hecke eigenform. Then,
by theorem (4.3.6) of \cite{St1}, under the assumption of the
proposition \ref{prop4},
\begin{eqnarray}
\Theta_{k,\chi}(\varphi_{f}^{-})=\theta_{k,\chi}^{\ast}(f)
\label{algebraic Shintani theta}
\end{eqnarray}
where $\varphi_f^-:=\frac{1}{\Omega_f^-}\cdot\psi_f^-$. Note that
$\theta^{\ast}_{k,\chi}$ depends on the choice of complex periods
and is defined only on Hecke eigenforms.

For $\hat{\kappa}\in B_{h,N}^{classical}(K)$, we similarly define
an $R_{\hat{\kappa}}$-linear map
\begin{eqnarray}
\phi_{\hat{\kappa}}: \cD_{B_{h,N}}=\cD\otimes_{A(\cX)}A(B_{h,N})
\to L_{k,\chi}(R_{\hat{\kappa}})\label{spB}
\end{eqnarray}
by replacing $\kappa$ by $\hat{\kappa}$ in (\ref{specialize}),
where $R_{\hat{\kappa}}:=\hat{\kappa}(A(B_{h,N}))$. Then the
following theorem due to G. Stevens (see [St3]) is one of the
major steps for the existence of an overconvergent Hecke-eigen
symbol.
\begin{thm}
Assume $\hat{\kappa} \in B_{h,N}^{classical}(K)$ of signature
$(k,\chi)$ and $h < k+1$. Then
\begin{quote}

(1) there is a Hecke-equivariant isomorphism

\begin{equation}
 H^1_c(\Gamma_{0}(Np), \cD_{B_{h,N}}
\otimes_{A(B_{h,N}),\hat{\kappa}} K)^{(\leq h)} \ {}^{\sim \atop
\to} \
 H^1_c(\Gamma_{0}(Np^m)), L_{k,\chi}(K))^{(\leq
 h)}\label{comparison1}
\end{equation}\\
induced from (\ref{spB}), where $m$ is the smallest positive
integer for which $\chi$ is defined modulo $Np^m$ and
$\otimes_{A(B_{h,N}),\hat{\kappa}}$ is the tensor product taken
with respect to $\hat{\kappa}:A(B_{h,N})\to K$.
\end{quote}

\begin{quote}
(2) There is a canonical identification

\begin{equation}
H^1_c(\Gamma_{0}(Np),
\cD_{B_{h,N}})\otimes_{A(B_{h,N}),\hat{\kappa}} K \simeq
H^1_c(\Gamma_{0}(Np), \cD_{B_{h,N}}
\otimes_{A(B_{h,N}),\hat{\kappa}} K)\label{comparison2}
\end{equation}\\
in the case $k >0$.
\end{quote}\label{thm4}
\end{thm}

In the remaining of this section, we will use the following
abbreviated notations:
\begin{eqnarray*}
H_c^1(\cD)^{\ast}&:=&H_{c}^{1}(\Gamma_{0}(Np),\cD_{B_{h,N}})^{(\leq
h),-}\\
H_c^1(L_{\hat{\kappa}})^{\ast}&:=&H_c^1(\Gamma_{0}(Np^m),
L_{k,\chi}(R_{\hat{\kappa}}))^{(\leq h),-}\\
\cR&:=&\cR_{h,N}
\end{eqnarray*}
where $\hat{\kappa}\in B^{classical}_{h,N}(K)$ is a $K$-point of
signature $(k,\chi)$, $0 \leq h < k+1$ and $m$ is the smallest
positive integer for which $\chi$ is defined modulo $Np^m$. Let
$m_{\hat{\kappa}}$ be a maximal ideal $\ker(\hat{\kappa})\subseteq
A(\BN)$. Robert Coleman proved that there is a point $\kappa \in
\ONc(\bC_p)$ which is unramified over the point $\hat{\kappa}\in
B_{h,N}^{classical}(K)$ in [Col1]. Put
$m_{\kappa}=\ker(\kappa)\subseteq \cR$.

\begin{cor}
There is a Hecke-equivariant $\cR/m_{\hat{\kappa}}\cR$-module
isomorphism between
\begin{eqnarray}
H_c^1(\cD)^{\ast}/m_{\hat{\kappa}}H_c^1(\cD)^{\ast} \simeq
H_c^1(L_{\hat{\kappa}})^{\ast}
\end{eqnarray}
induced from the map (\ref{spB}). The right hand side is an
$\cR(\hat{\kappa})$-module where
$\cR(\hat{\kappa}):=Im\big(R_{\hat{\kappa}}[T_n:n\in \bN] \to
End_K(H_c^1(L_{\hat{\kappa}})^{\ast})\big)$. We have a $K$-algebra
isomorphism from $\cR/m_{\hat{\kappa}}\cR$ to $\cR(\hat{\kappa})$
so that we can regard the right hand side a
$\cR/m_{\hat{\kappa}}\cR$-module by this $K$-algebra
isomorphism.\label{cor2}
\end{cor}
\begin{proof}
It follows from
\begin{eqnarray*}
H_c^1(\cD)^{\ast}/m_{\hat{\kappa}}H_c^1(\cD)^{\ast} & \simeq &
H_c^1(\cD)^{\ast}\otimes_{A(B_{h,N}),\hat{\kappa}}K\\
& \simeq & H_c^1(L_{\hat{\kappa}})^{\ast} \ (\textrm{by
(\ref{comparison1}) and (\ref{comparison2})}).
\end{eqnarray*}
\end{proof}

\begin{lem}
If $\kappa \in \ONc(\bC_p)$ is unramified over $\hat{\kappa}\in
B^{classical}_{h,N}(K)$, then $H_c^1(\cD)^{\ast}_{(\kappa)}$ is a
free $\cR_{(\kappa)}$-module of rank 1, where we use the subscript
$(\kappa_i)$ to denote the localization at $\kappa_i$ (using the
maximal ideal $m_{\kappa_i}$).\label{lem1}
\end{lem}


\begin{proof}

Note that $\cR/m_{\hat{\kappa}}\cR$ is a finite free $K$-algebra
and so it is the product of local $K$-algebras
$\prod_{i=1}^{d}\overline{\cR_{i}}$ ($d$ is some positive integer)
where each $\overline{\cR}_{i}$ is isomorphic to
$(\cR/m_{\hat{\kappa}}\cR)_{(\kappa_{i})}$=localization at
$\kappa_i$ of $\cR/m_{\hat{\kappa}}\cR$, where $\kappa_{i}'s$ are
points lying over $\hat{\kappa}$. Then there exists some $i>0$
such that $\kappa=\kappa_i$ which is unramified over
$\hat{\kappa}$. Let $m_{\kappa_i}$ be $\ker(\kappa_i)\subseteq
\cR$:
\begin{displaymath}
\xymatrix{m_{\kappa_i} \ar@{^{(}->}[r] & \cR  \ar[r]^{\kappa_i} & \overline{\bQ}_p \\
m_{\hat{\kappa}} \ar@{^{(}->}[r] & A(B_{h,N})  \ar[u]
\ar[ur]_{\hat{\kappa}}}.
\end{displaymath}
We denote the localization at $\kappa_i$ (using the maximal ideal
$m_{\kappa_i}$) by the subscript $(\kappa_i)$. If $N$ is
square-free, then it is known that
$H_c^1(L_{\hat{\kappa}})^{\ast}$ is a free
$\cR(\hat{\kappa})$-module of rank 1 by a slight extension of
Eichler-Shimura theory. If $N$ is not square-free, we still know
that $H_c^1(L_{\hat{\kappa}})^{\ast}_{(\kappa_{i})}$ is a free
$\overline{\cR_i}$-module of rank 1. So we conclude that
$\big(H_c^1(\cD)^{\ast}/m_{\hat{\kappa}}H_c^1(\cD)\big)^{\ast}_{(\kappa_i)}$
is a free $(\cR/m_{\hat{\kappa}}\cR)_{(\kappa_i)}$-module of rank
1 by the above corollary \ref{cor2}.

Then we have

\begin{eqnarray*}
\big(H_c^1(\cD)^{\ast}/m_{\hat{\kappa}}H_c^1(\cD)^{\ast}\big)_{(\kappa_i)}&\simeq&
H_c^1(\cD)^{\ast}_{(\kappa_i)}/m_{\hat{\kappa}}H_c^1(\cD)^{\ast}_{(\kappa_i)}\\
&\simeq&
H_c^1(\cD)^{\ast}_{(\kappa_i)}/m_{\kappa_{i}}H_c^1(\cD)^{\ast}_{(\kappa_i)}
\ (\textrm{because of unramifiedness}).
\end{eqnarray*}

So we have that
$H_c^1(\cD)^{\ast}_{(\kappa_i)}/m_{\kappa_{i}}H_c^1(\cD)^{\ast}_{(\kappa_i)}$
is a free $(\cR/m_{\hat{\kappa}}\cR)_{(\kappa_{i})}$-module of
rank 1. Hence it is a free
$\cR_{(\kappa_i)}/m_{\kappa_i}\cR_{(\kappa_i)}$-module of rank 1,
since $\kappa_{i}$ is unramified over $\hat{\kappa}$. Now if we
apply the Nakayama Lemma, we conclude that
$H_c^1(\cD)^{\ast}_{(\kappa_i)}$ is a free
$\cR_{(\kappa_i)}$-module of rank 1.
\end{proof}

\begin{rmk}
$H^1_c(\cD)^{\ast}$ doesn't have to be a free $\cR$-module of rank
1. Local free rank 1 property for every localization at prime
ideals doesn't imply the global rank 1 property in
general.\label{rem3}
\end{rmk}

\begin{thm}
Let $\kappa_{0} \in \Omega^{classical}_{h,N}(K)$ be any classical
$K$-point which is unramified over $B_{h,N}$. There is an
overconvergent Hecke eigensymbol $\Phi \in
H^1_c(\Gamma_0(Np),\cD_{B_{h,N}})^{(\leq h)}$ and a choice of
periods $\Omega_{\kappa} \in R_\kappa$ for $\kappa \in
\Omega^{classical}_{h,N}(K)$ such that \\

\begin{quote}
(1) $\Omega_{\kappa_0} \neq\ 0$
\end{quote}

\begin{quote}

(2) $\Phi_{\kappa}=\Omega_{\kappa}\varphi_{\kappa}^-$ for any
point $\kappa \in \Omega^{classical}_{h,N}(K)$.
\end{quote}\label{Hecke eigensymbol exists}

\end{thm}

\begin{proof}
\

Having Lemma \ref{lem1} in our hands, the theorem follows from the
same argument in the proof of the theorem 5.5 in \cite{St1}.
\end{proof}

\subsection{$p$-adic family of overconvergent half-integral weight
Hecke eigenforms via local overconvergent Shintani lifting}\

Let $\sigma:\Lambda_N \to \Lambda_N$ be the ring homomorphism
associated to the group homomorphism $t \to t^2$ on
$\bZ^{\times}_{p,N}$. Then the functoriality of the P. Berthelot's
construction of $K$-rigid analytic variety associated to
$\Lambda_N$ gives us the map $\sigma^{\ast}:\cX_N \to \cX_N$ and
consequently we get the ring homomorphism (using the same notation
$\sigma$)
\begin{eqnarray}
\sigma:A(\cX_N) \to A(X_N).
\end{eqnarray}
which commutes with the natural map $\Lambda_N \to A(X_N)$ and
$\sigma:\Lambda_N \to \Lambda_N$. If we use the $p$-adic Fourier
transform of Amice (\ref{Fourier}), the above $\sigma$ is in fact
same as $\sigma:\NS \to \NS$ in (\ref{double}).

\begin{defn}
We define
\begin{eqnarray}
\tilde{A}_{h,N}:=A(B_{h,N}) \otimes_{A(\cX_N),\sigma}A(\cX_N)
\end{eqnarray}
where $\otimes_{A(\cX_N),\sigma_N}$ is the completed tensor
product taken with respect to $\sigma:A(\cX_N) \to A(\cX_N)$ which
is the ring homomorphism given by the base change from
$\sigma:\Lambda_N \to \Lambda_N$.\label{def14}

\end{defn}

\begin{defn}
We define the overconvergent metaplectic $p$-adic Hecke algebra of
tame conductor $N$ with slope $\leq h$ by
\begin{equation}
\tilde{\cR}_{h,N}:=\cR_{h,N}
\otimes_{A(B_{h,N}),\sigma_{h}}\tilde{A}_{h,N}
\end{equation}
 where
$\otimes_{A(B_{h,N}),\sigma_{h}}$ is the completed tensor product
taken with respect to $\sigma_h:A(B_{h,N}) \to \tilde{A}_{h,N} $
which is the ring homomorphism given by the base change from
$\sigma:A(\cX_N) \to A(\cX_N)$.\label{def15}
\end{defn}
Now we can formulate the following diagram:
\begin{displaymath}
\xymatrix{ \tilde{\cR}_{h,N} & \cR_{h,N} \ar[l] & \cR_{h} \ar[l] \\
\tilde{A}_{h,N} \ar[u] & A(B_{h,N}) \ar[l]_{\sigma_{h}} \ar[u]  & A(B_{h}) \ar[l] \ar[u] \\
A(\cX_{N}) \ar[u] & A(\cX_{N}) \ar[l]_{\sigma} \ar[u] & A(\cX) \ar[u] \ar[l] \\
\Lambda_N \ar[u] & \Lambda_N \ar[l]_{\sigma} \ar[u] & \Lambda
\ar[u] \ar[l] }
\end{displaymath}

\begin{prop}
$\tilde{\cR}_{h,N}$ is a $K$-affinoid algebra.\label{prop5}
\end{prop}
\begin{proof}
The category of $K$-affinoid algebras is closed under completed
tensor product with respect to any contractive $K$-algebra
homomorphism (Proposition 10 on the page 225 of [BGR]). So the
proposition follows, since $\sigma_h$ is a contractive
homomorphism.
\end{proof}

We regard $\tilde{\cR}_{h,N}$ as a $\Lambda_N$-algebra by
equipping it with the structure homomorphism $\Lambda_N \to
\tilde{\cR}_{h,N}$ given by $f \mapsto
1\otimes_{A(B_{h,N}),\sigma_{h}}f$. Then we can easily check that
$\tilde{\cR}_{h,N}:=\cR_{h,N}
\otimes_{A(B_{h,N}),\sigma_{h}}\tilde{A}_{h,N}$ is isomorphic as
an $A(\cX_{N})$-algebra to
$\cR_{h,N}\otimes_{A(\cX_{N}),\sigma}A(\cX_{N})$ whose
$A(\cX_{N})$-algebra structure is given by $A(\cX_{N}) \to
\tilde{\cR}_N,\ \lambda \mapsto
1\otimes_{A(\cX_{N}),\sigma}\lambda$. So by the $p$-adic Fourier
transform (\ref{Fourier}), we have
\begin{eqnarray}
\tilde{\cR}_{h,N}\cong\cR_{h,N}\TNS \NS \label{Fouriertwist}
\end{eqnarray}
as a $\NS$-algebra. Note that the ring homomorphism
\begin{equation}
\cR_{h,N} \to \tilde{\cR}_{h,N}
\end{equation}
given by  $\alpha \mapsto \alpha \otimes_{A(\cX_{N}),\sigma} 1$ is
not a homomorphism of $A(\cX_{N})$-algebras. This is reflected in
the fact that the $K$-rigid analytic map induced by pullback on
the corresponding $K$-affinoid spaces
\begin{equation}
\tilde{\Omega}_{h,N} \to \Omega_{h,N}\label{pullback}
\end{equation}
does not preserve the signatures of classical points. Indeed, if
$\tilde{\kappa} \in \tilde{\Omega}_{h,N}^{classical}(K)$ has
signature $(k, \chi)$ and  lies over $\kappa \in
\Omega_{h,N}^{classical}(K)$, then the signature of $\kappa$ is
$(2k, \chi^2)$. If we let $\Omega_h$ (respectively,
$\tilde{B}_{h,N}$) be the $K$-affinoid space associated to $\cR_h$
(respectively, $\tilde{A}_{h,N}$), then we have the following
commutative diagram by pulling back the previous diagram:

\begin{displaymath}
\xymatrix{ \tilde{\Omega}_{h,N} \ar[r] \ar[d] & \Omega_{h,N} \ar[d] \ar[r] & \Omega_{h} \ar[d]  \\
\tilde{B}_{h,N} \ar[r]^{\sigma^*_{h}} \ar[d] & B_{h,N} \ar[d] \ar[r]  &  B_{h} \ar[d]\\
\cX_N \ar[r]^{\sigma^*} & \cX_{N}  \ar[r]  & \cX    }.
\end{displaymath}\\

Now we prepare for the final notations to define a local
overconvergent Shintani lifting over $\RN$. Let
\begin{eqnarray*}
\cD_{\cR_{h,N}}:=\cD_{\BN} \TAR
\end{eqnarray*}
where $S_0(Np)$ acts on $\cD_{\cR_{h,N}}$ through the first
factor. If $\kappa \in \Omega_{h,N}(K)$ is any $K$-point, then let
$R_\kappa:=\kappa(\cR_{h,N})$. For $\kappa \in \ONc(K)$ of
signature $(k,\chi)$, we define an $R_\kappa$-linear map
\begin{eqnarray}
\phi_\kappa: \cD_{\cR_{h,N}}=\cD\otimes_{A(\cX)}\cR_{h,N} \to
L_{k,\chi}(R_\kappa) \label{Rsp}
\end{eqnarray}
by
\begin{eqnarray*}
\phi_\kappa(\mu \otimes
r):=\kappa(r)\cdot\int_{\ZZ}\chi_p(x)\frac{(xY-yX)^k}{k!}d\mu(x,y)
\end{eqnarray*}
for $\mu \in \cD$ and $r \in \cR_{h,N}$, where we factor
$\chi=\chi_N \cdot \chi_p$ with $\chi_N$ defined modulo $N$ and
$\chi_p$ defined modulo a power of $p$. A simple computation
confirms that if $\chi$ is defined modulo $Np^m$ then
$\phi_\kappa$ is a $\Gamma_0(Np^m)$-module homomorphism, hence
induces a homomorphism on the cohomology groups:
\begin{eqnarray}
\phi_{\kappa,\ast}:\HR \to
H^1_c(\Gamma_0(Np^m),L_{k,\chi}(R_\kappa))
\end{eqnarray}which is Hecke-equivariant.

We state the following Lemma which is almost identical to the
proposition \ref{JQ} which is used crucially to construct a local
overconvergent Shintani lifting.
\begin{lem}
For each $Q \in \mathcal{F}_{Np}$ there is a unique
$\cR_{h,N}$-module homomorphism
\begin{eqnarray}
\tilde{J}^B_Q: \cD_{B_{h,N}} \otimes_{A(B_{h,N})} \cR_{h,N} \to
\tilde{\cR}_{h,N}
\end{eqnarray}\\
such that for all $\tilde{\kappa} \in
\tilde{\Omega}_{h,N}^{classical}(K)$ of signature $(k,\chi)$ lying
over $\kappa \in \Omega_{h,N}^{classical}(K)$, we have
\begin{eqnarray}
\tilde{\kappa}(\tilde{J}^B_{Q}(\mu))=\chi(Q) \cdot \langle
\phi_{\kappa}(\mu), Q^k \rangle
\end{eqnarray}\\
for $ \mu \in \cD_{B_{h,N}} \otimes_{A(B_{h,N})}
\cR_{h,N}=\cD_{\cR_{h,N}}$ and $\phi_{\kappa}(\mu)$ was defined in
\ref{Rsp}.\label{lem2}
\end{lem}

\begin{proof}
Once we define the appropriate map, the proofs of uniqueness and
interpolation properties are exactly same as the proposition
\ref{JQ}. So we concentrate on constructing the appropriate
$\RN$-linear map. We recall that there is a topological
$K$-algebra isomorphism (as $K$-Fr\'{e}chet space)
\begin{eqnarray*}
F:\cD(\bZ_p^{\times}) \simeq A(\cX)
\end{eqnarray*}
by the $p$-adic Fourier transform of Amice([Ami1, Ami2]).
 So using this isomorphism we get a well-defined $K$-linear map $\cD \to
 A(\cX)$ by $\nu \mapsto F(J_Q(\nu))$ where $J_Q$ is defined in (\ref{JQmap}).

 We let $[Q]_N:=[a]_N \in \Delta_N:=(\bZ/N)^{\times}$ for each
 quadratic form $Q(X,Y)=aX^2+bXY+cY^2 \in \mathcal{F}_{Np}$. Now
 we define $\cR_{h,N}$-homomorphism $J_{Q}: \cD_{B_{h,N}} \otimes_{A(B_{h,N})} \cR_{h,N} \to
\tilde{\cR}_{h,N}:=\cR_{h,N} \otimes_{A(B_{h,N}),\sigma_h}
\tilde{A}_{h,N}$ (remember that $
\tilde{A}_{h,N}:=A(B_{h,N})\otimes_{A(\cX_N),\sigma}A(\cX_N)$) by
the following:
\begin{displaymath}
\xymatrix{ \cD_{B_{h,N}} \otimes_{A(B_{h,N})} \cR_{h,N} \ar[r]^{\
\ \ \ \ \ \ \ \
\tilde{J}^B_Q} \ar@{=}[d] & \tilde{\cR}_{h,N} \ar@{=}[d]   \\
(\cD \otimes_{A(\cX)}A(B_{h,N}))\otimes_{A(B_{h,N})}\cR_{h,N}
\ar[r] & \cR_{h,N}
\otimes_{A(B_{h,N}),\sigma_h}(A(B_{h,N})\otimes_{A(\cX_N),\sigma}A(\cX_N))  \\
(\nu\otimes_{A(\cX)}f) \otimes_{A(B_{h,N})}r \ar[r] & r
\otimes_{A(B_{h,N}),\sigma_h}\Big(\sigma_h(f)\cdot \big(
1\otimes_{A(\cX_N),\sigma}[Q]_N\cdot
F(j_{Q}(\nu))\big)\Big)\ar@{|}[l] }
\end{displaymath}
for all $\nu \in \cD, f \in A(B_{h,N})$ and $r \in \cR_{h,N}$,
where $[Q]_N$ and $F(j_{Q}(\nu))$ are regarded as elements of
$A(\cX_N)$ by the natural inclusion $\Lambda_N \to A(\cX_N)$ and
$A(\cX) \to A(\cX_N)$ respectively and $\sigma_h(f)$ is the image
of
$f$ under $\sigma_h:A(B_{h,N}) \to \tilde{A}_{h,N}$.\\
Let $[t]$ be the image of $t$ under the natural map $\Lambda \to
A(\cX)$ and recall $A(\cX_N)$ acts on $\cD$. Since we have
$F(J_Q([t]\cdot\nu))=[t]^2\cdot F(J_Q({\nu}))$ and
$\sigma([t]\cdot f)=[t]^2\cdot \sigma(f)$ for all $t\in
\bZ_p^{\times}$, the map $\cD_{B_{h,N}} \to \tilde{A}_{h,N}$ by
$\nu\otimes_{A(\cX)}f \mapsto \sigma_h(f)\cdot
(1\otimes_{A(\cX_N),\sigma}[Q]_N\cdot F(j_{Q}(\nu)))$ is
$\bZ_p$-module homomorphism. Then $\tilde{J}^B_Q$ is actually an
$\cR_{h,N}$-linear extension of this map. If we identify $(\cD
\otimes_{A(\cX)}A(B_{h,N}))\otimes_{A(B_{h,N})}\cR_{h,N}$ with
$\cD \otimes_{A(\cX)} \cR_{h,N}$, then the map is given by
\begin{eqnarray}
\nu\otimes_{A(\cX)}r \mapsto r \otimes_{A(B_{h,N}),\sigma_h}
(1\otimes_{A(\cX_N),\sigma}[Q]_N\cdot F(j_{Q}(\nu))).
\end{eqnarray}
The interpolation property at classical points follows from the
same computation as in the proposition \ref{JQ}.
\end{proof}
We can define a local overconvergent Shintani lifting over $\BN$.

\begin{defn}
For each $\Phi \in H^1_{c}(\Gamma_{0}(Np), \cD_{\cR_{h,N}})$ and
each $Q \in \mathcal{F}$ we define
\begin{eqnarray}
J_B(\Phi, Q):=\tilde{J}^B_Q(\Phi(D_Q)) \in \tilde{\cR}_{h,N}
\end{eqnarray}
where $\Phi(D_Q)$ is the value of the cocycle representing $\Phi$
on $D_Q$. Note that this definition of $J(\Phi,Q)$ does not depend
on the representative cocycle, and depends only on the
$\Gamma_{0}(Np)$-equivalence class of $Q$.\label{defJQR-version}
\end{defn}

\begin{defn}
We define a local overconvergent Shintani lifting over $\BN$ to be
a map $\Theta_B:H^1_{c}(\Gamma_{0}(Np), \cD_{\cR_{h,N}}) \to
\tilde{\cR}_{h,N}[[q]]$ given by
\begin{displaymath}
\Theta_B(\Phi):=\left\{ \begin{array}{ll} \sum_{Q \in
\mathcal{F}_{Np}/\Gamma_{0}(Np)}J_B(\Phi,Q)q^{\delta_{Q}/Np} & \textrm{if $Np$ is odd,}\\
\\
\sum_{Q \in
\mathcal{F}_{Np}/\Gamma_{0}(Np)}J_B(\Phi,Q)q^{\delta_{Q}/4Np} &
\textrm{if $Np$ is even.}
\end{array} \right.\label{maindef2}
\end{displaymath}

\end{defn}

\begin{defn}
We say that $\Theta \in \RNT[[q]]$ is a universal overconvergent
half-integral weight modular form over $\BN$ if
\begin{eqnarray*}
\Theta=\Theta_B(\Phi)
\end{eqnarray*}
for some $\Phi \in \HBhN^{(\leq h)}$.\label{univmodB}
\end{defn}
Let me also give a definition of an overconvergent half-integral
weight modular form of weight $\kappa\in
\textrm{Hom}_{cts}(A(\BN),K)$.
\begin{defn}
We say that $\theta \in \bC_p[[q]]$ is an overconvergent
half-integral weight modular form of weight $\kappa\in
\textrm{Hom}_{cts}(A(\BN),K)$ if it can be written as
\begin{eqnarray*}
\theta=\Theta(\kappa)=\sum_{n=1}^{\infty}\big(r_n(\kappa)\TNS\alpha_n(\kappa)\big)q^n
\end{eqnarray*}
where $\Theta=\sum_{n=1}^{\infty}\big(r_n\TNS\alpha_n\big) q^n \in
\RNT[[q]]$ is a universal overconvergent half-integral weight
modular form.\label{overmodB}
\end{defn}

Motivated by the theorem \ref{main1}, it's tempting to define the
actions of $T_{l}$ for any rational prime $l$ and $T_{l,l}$ for
$l\nmid Np$ on a universal overconvergent half-integral weight
modular form over $\BN$
\begin{eqnarray*}
\Theta=\sum_{n=1}^{\infty}\big(r(n)\TNS b(n)\big)q^n \in \RNT[[q]]
\end{eqnarray*}
as follows:
\begin{eqnarray}
\Theta\big|T_l=
 \sum_{n=1}^{\infty}\Big(r(n) \TNS\big(b(nl^2)+(\frac{Np\cdot
 n}{l})[l]_N\delta_l \ast b(n)+l[l^2]_N\delta_{l^2}\ast
 b(\frac{n}{l^2})\big)\Big)q^n
\end{eqnarray}
\begin{eqnarray}
\Theta\big|T_{l,l}=\sum_{n=1}^{\infty}\Big(r(n) \TNS [l^2]_N \cdot
\delta_{l^2}\ast
 b(n)\Big)q^n.
\end{eqnarray}
Henceforth we can view $\RNT[[q]]$ as an $\cH$-module where $\cH$
was given in the Definition \ref{abstract Hecke algebra}. We warn
you that $T_n$-action on $\RNT[[q]]$ is NOT given by the
$\RN$-module structure of $\RNT$. Also notice that $\NS$ act on
$\RNT$ through the first factor. Then we can conclude that our
local overconvergent Shintani lifting $\Theta_B$ is $\cH$-linear
map by Theorem \ref{main1} except for operators containing $T_2$.
But if $N$ is even, then
$\Theta_B(\Phi\big|T_2)=\Theta_B(\Phi)\big|T_2$ will still be true
so that $\Theta$ is actually $\cH$-linear. Therefore we get the
following local version of the theorem \ref{main1}.
\begin{thm}
The map $\Theta_B$ is an $\cH$-module homomorphism except for
operators containing $T_2$.\label{h equivariance B}
\end{thm}
We can talk about a universal overconvergent half-integral weight
Hecke eigenform over $\BN$, since we defined Hecke action on
$\RNT[[q]]$. Theorem \ref{h equivariance B} implies that the map
$\Theta_B$ sends a Hecke eigensymbol to a universal overconvergent
half-integral weight Hecke eigenform over $\BN$. Now we focus on
the interpolation property of $\Theta_B$ and will find the
universal overconvergent half-integral weight Hecke eigen form
over $\BN$ corresponding to the Hecke eigensymbol we constructed
in the theorem \ref{Hecke eigensymbol exists}.
\begin{thm}
For each $\Phi \in H^1_{c}(\Gamma_{0}(Np), \cD_{\cR_{h,N}})$ and
for each classical point $\tilde{\kappa} \in
\tilde{\Omega}_{h,N}^{classical}(K)$ of signature $(k,\chi)$ lying
over $\kappa \in \Omega_{h,N}^{classical}(K)$ we have\\
\begin{eqnarray}
\Theta_B(\Phi)(\tilde{\kappa})\mid T_p^{m-1}
=\Theta_{k,\chi}(\Phi_{\kappa})
\end{eqnarray}\\
where $\Phi_{\kappa}=\phi_{\kappa,*}(\Phi)$ and $m$ is the
smallest positive integer for which $\chi$ is defined modulo
$Np^m$. Note that here we understand $T_p$ as the action on formal
$q$-expansion given by (\ref{Hecke p}).\label{interpolation of
Theta B}
\end{thm}

\begin{proof}
For $Np$ odd, this follows from the following calculation:
\begin{eqnarray*}
\Theta_B(\Phi)(\tilde{\kappa})\mid T_p^{m-1} &=& \Big( \sum_{Q\in
\mathcal{F}_{Np}/\Gamma_0(Np)}
\tilde{\kappa}(J(\Phi,Q))q^{\delta_Q/Np} \Big)|T_p^{m-1} \\
&=&\sum_{\substack{Q\in\mathcal{F}_{Np}/\Gamma_0(Np) \\
p^m|\delta_{Q}}}\chi(Q)\langle\phi_\kappa(\Phi(D_Q)),Q^k \rangle
q^{\delta_Q/Np^m}\\
&=&\sum_{Q\in\mathcal{F}_{Np^m}/\Gamma_0(Np^m)}\chi(Q)\langle\Phi_\kappa(D_Q),Q^k
\rangle q^{\delta_Q/Np^m} \\
&=&\sum_{Q\in\mathcal{F}_{Np^m}/\Gamma_0(Np^m)}J_{k,\chi}(\Phi_\kappa,Q
) q^{\delta_Q/Np^m} \\
&=& \Theta_{k,\chi}(\Phi_k).
\end{eqnarray*}
The same computation applies to the $Np$ even case.
\end{proof}

We state the main theorem on the existence of the universal
overconvergent half-integral weight Hecke eigenform over $\BN$ and
prove it.

\begin{thm}

Let $\kappa_0 \in \Omega_{h,N}^{classical}(K)$ be a fixed
classical point unramified over $B_{h,N}$. Then there is a formal
$q$-expansion $\Theta =\sum_{n=1}^{\infty}\beta_n q^n \in
\tilde{R}_{h,N}[[q]]$ and a choice of periods $\Omega_{\kappa} \in
K$, for $\kappa \in \Omega_{h,N}^{classical}(K)$, with the
following properties:

\begin{quote}
(1) $\Omega_{\kappa_{0}}\neq 0$
\end{quote}
\begin{quote}
(2) For every classical point $\tilde{\kappa} \in
\tilde{\Omega}_{h,N}^{classical}(K)$ of signature $(k,\chi)$,\\
\begin{eqnarray*}
\Theta(\tilde{\kappa}):=\sum_{n=1}^{\infty}\beta_n(\tilde{\kappa})q^n
\in S_{k+\frac{3}{2}}(\Gamma_{0}(4Np^{m}), \chi^{\ast}; K)
\end{eqnarray*}
\\
where $m$ is the smallest positive integer for which $\chi$ is
defined modulo $Np^m$ and $\chi^{\ast}$ is defined by
$\chi^{\ast}(d)=\chi(d)(\frac{(-1)^{(k+1)}Np}{d})$.\label{thm8}
\end{quote}

\begin{quote}
(3) If $\kappa$ is the image of $\tilde{\kappa} \in
\tilde{\Omega}_{h,N}^{classical}(K)$ of signature $(k,\chi)$ under
the map (\ref{pullback}), then \\
\begin{eqnarray}
\Theta(\tilde{\kappa})=\Omega_{\kappa} \cdot
\theta^{\ast}_{k,\chi} (\mathbf{f}_{h,N}(\kappa)) \mid T_p^{1-m}
\end{eqnarray}
\\
where $\theta^{\ast}_{k,\chi} (\mathbf{f}_{h,N}(\kappa)) \mid
T_p^{1-m}:=\alpha_p(\kappa)^{1-m} \cdot \theta^{\ast}_{k,\chi}
(\mathbf{f}_{h,N}(\kappa)) \mid T_p^{m-1}$.
\end{quote}
\end{thm}
\begin{proof}
There is an Hecke eigensymbol $\Phi \in \HBhN^{(\leq h)}$ and a
choice of periods $\Omega_\kappa \in R_\kappa$ by theorem
\ref{Hecke eigensymbol exists}. We will use $\Phi\otimes1$ to
prove the existence of a formal $q$-expansion $\Theta$. In order
for that, we should make sure that $\Phi\otimes1 \in \HR$. This
follows from the flatness of $\cR_{h,N}$ as $A(B_{h,N})$-module,
which was guaranteed by proposition \ref{prop3}:
\begin{displaymath}
\xymatrix{ \HBhN^{(\leq h)} \ar@{^{(}->}[r] & \HBhN}
\end{displaymath}
preserves the injectiveness after tensoring $\TAR$
\begin{displaymath}
\xymatrix{ \HBhN^{(\leq h)}\TAR \ar@{^{(}->}[r] & \HBhN\TAR \ar@{^{(}->}[d]\\
 & \HR}
\end{displaymath}
where the second map $\HBhN\TAR \hookrightarrow \HR$ is given by
sending $\phi\otimes_{A(B_{h,N})} r$ to $D \mapsto
\phi(D)\otimes_{A(B_{h,N})}r$ for $D\in \Delta_0$, which is easily
checked to be a well-defined injective homomorphism.\\
Now we can define
\begin{displaymath}
\Theta:=\Theta_B(\Phi\otimes1)=\sum_{n=1}^{\infty}\beta_nq^n\in\tilde{\cR}_{h,N}[[q]].
\end{displaymath}
Then by (1) of theorem \ref{Hecke eigensymbol exists}, (1)
follows. Let's prove (2) and (3). Let $\tilde{\kappa} \in
\tilde{\Omega}_{h,N}^{classical}(K)$ be a classical $K$-point of
signature $(k,\chi)$ and let $\kappa$ be its image in $\ONc(K)$
under the map (\ref{pullback}). By the theorem \ref{interpolation
of Theta B} and (\ref{algebraic Shintani theta}) we have
\begin{eqnarray*}
\Theta(\Phi)(\tilde{\kappa})\mid T_p^{m-1}
&=&\Theta_{k,\chi}(\Phi_{\kappa})\\
&=&\Theta_{k,\chi}(\Omega_{\kappa} \cdot \varphi_\kappa^-)\\
&=&\Omega_{\kappa}\cdot
\theta_{k,\chi}^{\ast}(\mathbf{f}_{h,N}(\kappa))
\end{eqnarray*}
By the Shintani's theorem \ref{Shintani's main thm}, for $m=1$, it
follows that $\Theta(\tilde{\kappa})|T_p^2=\Omega_{\kappa}\cdot
\theta_{k,\chi}^{\ast}(\mathbf{f}_{h,N}(\kappa)|T_p)=\Omega_{\kappa}\cdot
\theta_{k,\chi}^{\ast}(\alpha_p(\kappa)
\cdot\mathbf{f}_{h,N}(\kappa))=\alpha_p(\kappa)\cdot
\Theta(\tilde{\kappa})$. So we can conclude that
$(\beta_{np^2})(\tilde{\kappa})=\alpha_p(\kappa)\cdot\beta_{n}(\tilde{\kappa})$
for every $n$ and every $\tilde{\kappa} \in
\tilde{\Omega}_{h,N}^{classical}(K)$ of signature $(k,\chi)$ with
$\chi$ defined modulo $Np$. Hence we have that
\begin{eqnarray}
\beta_{np^2}=\alpha_p\cdot \beta_n \in \tilde{\cR}_{h,N}
\label{Thetarelation}
\end{eqnarray}
for all $n \geq 1$. Now we apply the Hecke operator $T_p^{m-1}$ to
$\Theta(\Phi)(\tilde{\kappa})\mid T_p^{m-1}=\Omega_{\kappa}\cdot
\theta_{k,\chi}^{\ast}(\mathbf{f}_{h,N}(\kappa))$ and multiply by
$\alpha_p(\kappa)^{1-m}$ to obtain
\begin{eqnarray*}
\Theta(\tilde{\kappa})&=&\alpha_p(\kappa)^{1-m}\cdot
\Theta(\tilde{\kappa})|T_p^{2(m-1)}\ (\textrm{by} \ (\ref{Thetarelation}))\\
&=& \alpha_p(\kappa)^{1-m}\cdot
\Omega_{\kappa}\cdot\theta_{k,\chi}^{\ast}\big(\mathbf{f}_{h,N}(\kappa)\big)|T_p^{m-1}\\
&=:&\Omega_{\kappa}\cdot\theta_{k,\chi}^{\ast}\big(\mathbf{f}_{h,N}(\kappa)\big)|T_p^{1-m}.
\end{eqnarray*}
This proves the part (3). The classical Shintani's theorem 2.3
tells us that
\begin{eqnarray*}
\theta^{\ast}_{k,\chi}(\mathbf{f}_{h,N}(\kappa))\in
S_{k+\frac{3}{2}}(\Gamma_{0}(4Np^{m}), \chi'; K)
\end{eqnarray*}
where $\chi'$ is the character of $(\bZ/4Np^m\bZ)^{\times}$
defined by $\chi'(d):=\chi(d)\cdot
\big(\frac{(-1)^{k+1}Np^m}{d}\big)$.(see (2.2)) If we apply $T_p$
to $\theta^{\ast}_{k,\chi}(\mathbf{f}_{h,N}(\kappa))$, then it
multiplies the Nebentype by $\big(\frac{p}{\cdot}\big)$ by
proposition 1.5 of \cite{Shm2}. Therefore we have that
\begin{eqnarray*}
\theta^{\ast}_{k,\chi}(\mathbf{f}_{h,N}(\kappa))|T_p^{1-m}\in
S_{k+\frac{3}{2}}(\Gamma_{0}(4Np^{m}), \chi^{\ast}; K).
\end{eqnarray*}
where $\chi^{\ast}(d)=\chi(d)\big(\frac{(-1)^{(k+1)}Np}{d}\big)$.

So we conclude that $\Theta(\tilde{\kappa}) \in
S_{k+\frac{3}{2}}(\Gamma_{0}(4Np^{m}), \chi^{\ast}; K)$. (part
(2)) This finishes the proof of the main theorem.
\end{proof}

\begin{rmk}
$\Omega_{h,N}$ (respectively $\tilde{\Omega}_{h,N}$) can be viewed
as a $K$-affinoid subvariety of the eigencurve (respectively
half-integral eigencurve). So our construction in fact gives us
the local $K$-rigid analytic map from integral eigencurve to
half-integral eigencurve.
\end{rmk}

We finish the paper summarizing the whole picture of the
overconvergent Shintani lifting and its interpolation properties
by the following commutative diagram:

\begin{displaymath}
\xymatrix{ \Hc \ar[r]^{\ \ \ \ \Theta} \ar[d] & \DN[[q]] \ar[d] \\
\HBhN^{(\leq h)} \ar[d]^{\phi_{\kappa,*}} \ar[r]^{\ \ \ \ \ \ \ \ \ \Theta_B} & \RNT[[q]] \ar[d]^{\tilde{\kappa}} \\
\HLs  \ar[d]  \ar[dr]^{\ \ \ \ \ \Theta_{k,\chi}} & \bC_p[[q]] \ar[d]^{T_p^{m-1}} \\
\HLsp^{\pm} \ar[r]_{\ \ \ \ \ \ \ \ \ \ \ \ \ \ \Theta_{k,\chi}} & \bC_p[[q]] \\
\Cusp \ar[u]^{\cong} \ar[r] \ar[r]^{\theta_{k,\chi}} & \Hcusp.
\ar[u]}
\end{displaymath}



\end{document}